\documentclass{amsart}
\usepackage{amssymb,amsmath,amscd}

    \def\nz{\hbox{\text{0}{ \raise 3pt \hbox{\kern -10pt
    \vrule width 8pt height 0.4 pt }\kern 2pt}}}

\begin{document}
\title{JACOBIAN PAIRS}
\thanks{
Mathematics Department, Purdue University, 
West Lafayette, IN 47907, USA,\\
\hspace*{5mm}e-mail: ram@cs.purdue.edu\\
\hspace*{4mm}
Universit\'e d'Angers, Math\'ematiques, 49045 Angers cedex 01, France;\\
\hspace*{4mm}
e-mail: assi@univ-angers.fr 
Abhyankar's work was partly supported by NSF Grant DMS 97-32592 and
NSA grant MDA 904-99-1-0019.
During the development of this work, Assi visited Purdue University; he would
like to thank that institution for hospitality and support.}
\author{By Shreeram S. Abhyankar and Abdallah Assi}
\date{}
\begin{abstract}
We study meromorphic jacobian pairs, i.e., pairs of 
polynomials in one variable, with coefficients meromorphic series in a second
variable, whose jacobian relative to the two variables depends only on
the second variable. We pose two meromorphic jacobian conjectures about such
pairs one of which is in terms of an invariant of the pair which we call the 
beta invariant. These conjectures are shown to imply the bivariate algebraic
jacobian conjecture which predicts that two bivariate polynomials
generate the polynomial ring if their jacobian is a nonzero constant.
Actually we define the beta invariant for any two meromorphic curves, i.e.,
polynomials in one variable with coefficients which are meromorphic series
in another variable. One of our basic techniques is the beta-jacobian identity
which relates the beta invariant to the jacobian of the meromorphic curves.
When the pair of meromorphic curves consists of a curve and its derivative,
the beta invariant is reduced to the beta-bar invariant of the curve, 
and the beta-jacobian identity is reduced to the betabar-derivative identity.
When the curve is analytic, i.e., given by a bivariate power series, the
betabar-derivative identity gives rise to the conductor-derivative identity
which is related to Dedekind's formula expressing the derivative ideal as
a product of the conductor and the different. In turn, when the curve is
algebraic, i.e., when the power series is a bivariate polynomial, the
betabar-derivative identity gives rise to the betabar-genus identity
which connects the betabar invariant of the curve to its genus.
In the complex case the betabar invariant is related to 
the rank of the first homology group and the numbers of Milnor and Tjurina.
As another technique for studying the jacobian conjectures we revisit
the Newton polygon.
\end{abstract}
\maketitle

\centerline{}

\centerline{\bf Section 1: Introduction}

\centerline{}

Let $F=F(X,Y)$ and $G=G(X,Y)$ be meromorphic curves, i.e.,
let $F$ and $G$ be members of $R=k((X))[Y]=$ the polynomial ring in
$Y$ over the meromorphic series field $k((X))$ where $k$ is an algebraically
closed field of characteristic zero.
We want to study the jacobian $J(F,G)=F_XG_Y-F_YG_X$ of $F$ and $G$
relative to $X$ and $Y$ where subscripts indicate partial derivatives.
To do this we shall consider factors of $F$ in
$R^{\natural}=$ the set of all irreducible monic
polynomials (of positive degree) in $Y$ over $k((X))$. 
In general we shall use the terminology of \cite{AA1} and \cite{AA2}.

Let
$$
N=\text{deg}_Y F
$$
and note that: 
if $F\ne 0$ then $N$ is a nonnegative integer, and if $F=0$ then 
$N=-\infty$.
Let $F_0=F_0(X)\in k((X))$ be defined by saying that if $F=0$ then $F_0=0$,
and if $F\ne 0$ then $F_0=$ the coefficient of $Y^N$ in $F$. We say that 
$F$ is {\bf $k$-monic} to mean that $0\ne F_0\in k$. We can write
$$
F=F_0\prod_{1\le j\le\chi(F)}F_j\quad\text{ where }\quad
F_j=F_j(X,Y)\in R^{\natural}\text{ with deg}_YF_j=N_j
$$ 
for $1\le j\le\chi(F)$, and $\chi(F)$ is a nonnegative integer with:
$\chi(F)>0\Leftrightarrow N>0$.
We call $\chi(F)$ the {\bf branch number} of $F$.
For any integer $\nu>0$ which is divisible by 
$N_1,\dots,N_{\chi(F)}$, by Newton's Theorem we can write
$$
F(X^{\nu},Y)=F_0(X^\nu)\prod_{1\le i\le N}(Y-z_i(X))\quad\text{ with }\quad
z_i(X)\in k((X))
$$ 
for $1\le i\le N$, and upon letting
$$
M=\text{deg}_Y G
$$
we define the {\bf intersection multiplicity} 
of $F$ with $G$ by putting
$$
\text{int}(F,G)=\begin{cases}
\frac{M}{\nu}\text{ord}_XF_0(X^{\nu})
+(1/\nu)\sum_{i=1}^{N}\text{ord}_X G(X^{\nu},z_i(X))
\text{ if }F\ne 0\ne G\\
0\;\quad\text{ if }F=0\ne G\in k((X))\text{ or }G=0\ne F\in k((X))\\
\infty\quad\text{ if }F=0\ne G\not\in k((X))\text{ or }G=0\ne F\not\in k((X))\\
\infty\quad\text{ if }F=0=G\\
\end{cases}
$$
and by noting that this is independent of the choice of $\nu$.
We follow the convention that a sum of a finite number of
quantities is $\infty$ if each of them is either a real number or $\infty$ 
and at least one of them is $\infty$. Moreover the product of $\infty$ and
a positive real number is $\infty$. Also the sum of an empty family is
zero, and the product of an empty family is one.
Note that if $FG\ne 0$ then
$$
\text{int}(F,G)=\text{ord}_X(\text{Res}_Y(F,G))
$$
where Res$_Y(F,G)$ is the $Y$-resultant of $F$ and $G$, and hence always
$$
\text{int}(F,G)=\text{int}(G,F)=\text{an integer or }\infty.
$$
We define GCD$(F,G)=1$ to mean that int$(F,G)\ne\infty$ and we note that then
$$
\text{GCD}(F,G)=1\Leftrightarrow
\begin{cases}
\text{either $F\ne 0\ne G$ and $G\not\in F_jR$ for $1\le j\le\chi(F)$,}\\
\text{or }F=0\ne G\in k((X)),\\
\text{or }G=0\ne F\in k((X)).\\
\end{cases}
$$

For a moment suppose that $F\in R^{\natural}$; then for $1\le i\le N$ we 
have int$(F,G)=\text{ord}_X G(X^{\nu},z_i(X))$ and:
int$(F,G)\ge 0\Leftrightarrow\text{ord}_X(G(X^{\nu},z_i(X))-\lambda)>0$
for some $\lambda\in k$; 
note that if this $\lambda$ exists then it is 
unique and is independent of $i$ as well as $\nu$, and we call it the 
{\bf residue} of $G$ at $F$ and denote it by res$(F,G)$; also note that the 
said $\lambda$, if it exists, 
can also be characterized as the unique element of 
$k$ such that int$(F,G-\lambda)>0$;
if int$(F,G)<0$ then we put res$(F,G)=\infty$. Thus 
$$
\begin{cases}
\text{if }F\in R^{\natural}\text{ then }
\lambda=\text{res}(F,G)\in k\cup\{\infty\}\text{ is such that:}\\
\lambda=\infty\Leftrightarrow\text{int}(F,G)<0,\text{ and}\\
\lambda\in k\;\Leftrightarrow\text{int}(F,G)\ge 0
\Leftrightarrow\text{int}(F,G-\lambda)>0.
\end{cases} 
$$

In the general case  of $F\in R$, for any $\lambda\in k$
we define the {\bf $\lambda$-alpha invariant} of $G$ relative to $F$ and
the {\bf $\lambda$-beta invariant} of $G$ relative to $F$ by
putting
$$
\alpha_{\lambda}(F,G)=\{j:1\le j\le\chi(F)\text{ with res}(F_j,G)=\lambda\}
$$ 
and
$$
\beta_{\lambda}(F,G)=\sum _{j\in\alpha_{\lambda}(F,G)}
\text{int}(F_j,G-\lambda)
$$ 
respectively, and we define the {\bf alpha invariant} of $G$ relative 
to $F$ and the {\bf beta invariant} of $G$ relative to $F$ by putting
$$
\alpha(F,G)=\{\lambda\in k:\alpha_{\lambda}(F,G)\ne\emptyset\}
$$
and
$$
\beta(F,G)=\sum_{0\ne\lambda\in\alpha(F,G)}\beta_{\lambda}(F,G)
$$
respectively.

Note that although the invariant $\alpha_{\lambda}(F,G)$ depends on the way
the factors $F_j$ of $F$ are labelled, the other three invariants do not.
We could have avoided this dependence by defining $\alpha_{\lambda}(F,G)$
to be the obvious effective divisor of $R$, where by an effective divisor
of $R$ we mean a nonnegative-integer-valued map of $R^{\natural}$ with
finite support. We would have then started by associating to $F$ the 
effective divisor $F_D$ which sends any $\Phi\in R^{\natural}$ to the
number of $j$ with $1\le j\le\chi(F)$ for which $F_j=\Phi$, and so on.

Note that for any $\lambda\in k$ we have
$$
|\alpha_{\lambda}|\le\chi(F)<\infty
$$ 
where $|\;|$ denote size (= cardinality) and
$$
\begin{cases}
\beta_{\lambda}(F,G)
\text{ is a nonnegative integer or $\infty$ such that:}\\
\beta_{\lambda}(F,G)=0\Leftrightarrow\alpha_{\lambda}(F,G)=\emptyset,
\text{ and }\\
\beta_{\lambda}(F,G)=\infty\Leftrightarrow G-\lambda\in F_jR
\text{ for some }j\in\alpha_{\lambda}(F,G).
\end{cases}
$$
Moreover
$$
\alpha(F,G)=\{\lambda\in k:\text{res}(F_j,G)=\lambda\text{ for some }j\}
$$
and hence
$$
|\alpha(F,G)|\le\chi(F)
-|\{j:1\le j\le\chi(F)\text{ with int}(F_j,G)<0\}|<\infty
$$
and
$$
\begin{cases}
\beta(F,G)\text{ is a nonnegative integer or $\infty$ such that:}\\
\beta(F,G)=0\Leftrightarrow\alpha(F,G)\subset\{0\},
\text{ and }\\
\beta(F,G)=\infty\Leftrightarrow G-\lambda\in F_jR
\text{ for some $0\ne\lambda\in\alpha(F,G)$ and }j\in\alpha_{\lambda}(F,G).
\end{cases}
$$
Defining the {\bf minimal-intersection} of $F$ with $G$ by putting
$$
\text{minint}(F,G)=\text{min}_{\mu\in k}\text{int}(F,G-\mu)
$$
(where min stands for glb = greatest lower bound) we see that
$$
\begin{cases}
\text{minint}(F,G)=\text{an integer or }\infty,\\
\text{and: minint}(F,G)=\infty\Leftrightarrow F=0\ne G\not\in k((X))
\end{cases}
$$
and
\begin{align*}
\alpha(F,G)&=\{\lambda\in k:\text{int}(F,G-\lambda)>\text{minint}(F,G)\}\\
&=\{\lambda\in k:\text{int}(F,G-\lambda)\ne\text{minint}(F,G)\}
\end{align*}
and
\begin{align*}
\beta(F,G)&=\sum_{0\ne\lambda\in\alpha(F,G)}
[\text{int}(F,G-\lambda)-\text{minint}(F,G)]\\
&=\text{a nonnegative integer or }\infty.
\end{align*}

By a {\bf relative irregular value} of $G$ (on $F$) we mean an element of
$\alpha(F,G)$.  We call $\alpha(F,G)$ the {\bf relative irregular value set} 
of $G$ (on $F$), and we call $\alpha_{\lambda}(F,G)$ the 
{\bf relative irregular $\lambda$-label set} of $G$ (on $F$). We call 
$\beta(F,G)$ the {\bf relative excess intersection} of $G$ (on $F$), and we
call $\beta_{\lambda}(F,G)$ the {\bf relative excess intersection} of 
$G$ (on $F$) at $\lambda$.
[We may think of the $\alpha_{\lambda}$'s, with $\lambda\ne 0$, 
as branches of $F$ close to $G$, each of them giving a residue which when 
subtracted from $G$ gives a higher intersection multiplicity, and
the $\beta$'s are the sums of these higher intersection multiplicities.]

For any $\lambda\in k$ we put 
$$
\overline\alpha_{\lambda}(F)=\alpha_{\lambda}(F_Y,F)
$$ 
and
$$
\overline\beta_{\lambda}(F)=\beta_{\lambda}(F_Y,F)
$$ 
and call these the {\bf $\lambda$-alphabar invariant} of $F$ and the
{\bf $\lambda$-betabar invariant} of $F$ respectively. Also we put
$$
\overline\alpha(F)=\alpha(F_Y,F)
$$ 
and 
$$
\overline\beta(F)=\beta(F_Y,F)
$$
and call these the {\bf alphabar invariant} of $F$ and the 
{\bf betabar invariant} of $F$ respectively. 
By an {\bf irregular value} of $F$ we mean an element of
$\overline\alpha(F)$. We call $\overline\alpha(F)$ the 
{\bf irregular value set} of $F$, and we call $\overline\alpha_{\lambda}(F)$ 
the {\bf irregular $\lambda$-label set} of $F$. We call $\overline\beta(F)$ 
the {\bf excess intersection} of $F$, and we call $\overline\beta_{\lambda}(F)$ 
the {\bf excess intersection} of $F$ at $\lambda$.
The objects $\overline\alpha(F)$ and $\overline\beta(F)$ were considered
in \cite{Ass} where they were denoted by $I(F)$ and $A_F$ respectively. 

In Section 2 we shall prove an identity for $\beta(F,G)$ involving
int$(F,J(F,G))$ which we call the {\bf beta-jacobian identity}, and from 
this we shall deduce an identity for $\beta(F,G)$ involving
int$(F,F_Y)$ which we call the {\bf beta-derivative identity}, and in turn
from this we shall deduce an identity for $\overline\beta(F)$ involving 
int$(F,F_Y)$ which we call the {\bf betabar-derivative identity} 
and which was directly proved in \cite{Ass}.

The beta-jacobian identity gives a geometric
characterization of the beta invariant
without digging into the branch structure.
The beta-derivative identity says that if $G$ has a special relationship
with $F$ then the  beta-invariant has a simpler expression than the
one given by the beta-jacobian identity.

In Sections 3 and 4 we shall deduce further corollaries of the beta-jacobian
identity for analytic curves (given by power series in $X,Y$) and
algebraic curves (given by polynomials in $X,Y$) respectively.
In Sections 5 and 6 we shall consider curves with no irregular value and
one irregular value respectively.
In Section 7 we shall make two meromorphic jacobian conjectures 
and indicate how they imply the usual
bivariate algebraic jacobian conjecture, which predicts that if $F$ and
$G$ are algebraic curves, i.e., if they are members of the bivariate
polynomial ring $k[X,Y]$, then:
$0\ne J(F,G)\in k\Leftrightarrow k[X,Y]=k[F,G]$; note that the implication
$\Leftarrow$ is obvious since it follows by the chain rule for jacobians.
In Section 8 we shall settle some cases of the first 
meromorphic jacobian conjecture. In Sections 9 and 10, as a tool for studying 
the second meromorphic jacobian conjecture,
we shall revisit the Newton Polygon.
In particular we shall prove a parallelness property for the Newton Polygons
of two meromorphic curves whose jacobian depends only on one of the variables.
In Section 11 we shall relate some of the invariants
studied in Sections 1 to 6 with homology rank and the numbers of Milnor
and Tjurina.

It is a pleasure to thank Avinash Sathaye for many stimulating 
conversations about the material of this paper.

\centerline{}

\centerline{\bf Section 2: The beta-jacobian identity} 

\centerline{}

Let $E$ be the square free part of $F\in R$, i.e.,
in the notation of Section 1,
$E=F_0\prod_{1\le i\le b(F)}F_{a(i)}$ where $a(1),\dots,a(b(F))$ are
distinct integers amongst $1,\dots,\chi(F)$ such that 
$F_{a(1)},\dots,F_{a(b(F))}$ are all the distinct elements amongst
$F_1,\dots,F_{\chi(F)}$; let us call $E$ the {\bf radical} of $F$ and denote
it by rad$(F)$. Recall that for $G\in R$ we have: 
GCD$(F,G)=1\Leftrightarrow\text{int}(F,G)\ne\infty$.
Now if $F\ne 0$ then clearly: GCD$(F,F_Y)=1\Leftrightarrow\text{rad}(F)=F$. 
Using these concepts, we shall now prove:

\centerline{}

{\bf The beta-jacobian identity (2.1).} {\it Let $F\in R$ be $k$-monic of
$Y$-degree $N$, and let $G\in R$ be such that {\rm GCD}$(F,G-c)=1$ for all 
$c\in k$. Then for $E={\rm rad}(F)$ we have
$$
{\rm int}(F,J(E,G))=\text{int}(F,G)+{\rm int}(F,E_Y)-N+\beta(F,G)
$$
where each term is an integer.}

\centerline{}

PROOF. For $N=0$ this is obvious because each term is reduced to $0$.
For a moment suppose $N=1$. Then $E(X,Y)=F(X,Y)=F_0(Y-z(X))$ with
$0\ne F_0\in k$ and $z(X)\in k((X))$. Clearly
$J(E,G)=-F_0z_X(X)G_Y(X,Y)-F_0G_X(X,Y)$ and by substituting $z(X)$
for $Y$ in the Right Hand Side of this equation it becomes equal to
$-F_0H_X(X)$ where $H(X)=G(X,z(X))$. Therefore 
int$(F,J(E,G))=\text{ord}_X H_X(X)$. Also obviously
int$(F,G)=\text{ord}_X H(X)$ and int$(F,E_Y)=0$. 
If ord$_XH(X)\ne 0$ then
$\beta(F,G)=0=\text{ord}_X H_X(X)-\text{ord}_X H(X)+1$,
and if ord$_XH(X)=0$ then $H(X)=\lambda+\overline H(X)$ where
$0\ne\lambda\in k$ and $0\ne \overline H(X)\in k[[X]]$ with
$\overline H(0)=$ and hence again
$\beta(F,G)=\text{ord}_X\overline H(X)
=\text{ord}_X H_X(X)-\text{ord}_X H(X)+1$. Thus always
$\beta(F,G)=\text{ord}_X H_X(X)-\text{ord}_X H(X)+1$ and hence
$\text{int}(F,J(E,G))=\text{int}(F,G)+\text{int}(F,E_Y)-N+\beta(F,G)$
where each term is an integer.

To prove the general case, let the notation be as in Section 1. 
Let $\widehat E=\widehat E(X,Y)=E(X^{\nu})$
and $\widehat F=\widehat F(X,Y)=F(X^{\nu})$
and $\widehat G=\widehat G(X,Y)=G(X^{\nu})$. 
Let $\widehat F_0=F_0$ and $\widehat N_0=0$, and for $1\le i\le N$
let $\widehat F_i=\widehat F_i(X,Y)=Y-z_i(X)$ and $\widehat N_i=1$.
Let $\{u(1,1),\dots,u(1,v(1))\},\dots,\{u(w,1),\dots,u(w,v(w))\}$ be
a partition of $\{0,\dots,N\}$ into disjoint nonempty subsets
such that $\widehat F_{u(1,1)},\dots,F_{u(w,1)}$ are exactly all the
distinct elements amongst $\widehat F_0,\dots,\widehat F_N$, and
for all $i,j,j'$ we have $\widehat F_{u(i,j)}=\widehat F_{u(i,j')}$.
Then
$$
\widehat F=\prod_{1\le i\le w}\prod_{1\le j\le v(i)}\widehat F_{u(i,j)}
\quad\text{ and }\quad
\widehat E=\prod_{1\le i\le w}\widehat F_{u(i,1)}.
$$
By the above two cases we see that for $1\le i\le w$ and $1\le j\le v(i)$
we have
\begin{align*}
\text{int}(\widehat F_{u(i,j)},J(\widehat F_{u(i,1)},\widehat G))
=&\text{int}(\widehat F_{u(i,j)},\widehat G)
+\text{int}(\widehat F_{u(i,j)},\widehat F_{u(i,1)Y})\\
&-\widehat N_{u(i,j)}+\beta(\widehat F_{u(i,j)},\widehat G)
\end{align*}
where each term is an integer. Adding 
int$(\widehat F_{u(i,j)},\widehat E/\widehat F_{u(i,1)})$ to both 
sides of the above
equation, and noting that by various product rules and such we have
\begin{align*}
&\text{int}(\widehat F_{u(i,j)},\widehat E/\widehat F_{u(i,1)})
+\text{int}(\widehat F_{u(i,j)},J(\widehat F_{u(i,1)},\widehat G))\\
&=\text{int}(\widehat F_{u(i,j)},
(\widehat E/\widehat F_{u(i,1)})J(\widehat F_{u(i,1)},\widehat G))\\
&=\text{int}
\left(\widehat F_{u(i,j)},\sum_{0\le i'\le w\text{ with }i'\ne i}
(\widehat E/\widehat F_{u(i,1)})J(\widehat F_{u(i',1)},\widehat G)\right)\\
&=\text{int}(\widehat F_{u(i,j)},J(\widehat E,\widehat G))
\end{align*}
and
\begin{align*}
&\text{int}(\widehat F_{u(i,j)},\widehat E/\widehat F_{u(i,1)})
+\text{int}(\widehat F_{u(i,j)},\widehat F_{u(i,1)Y})\\
&=\text{int}(\widehat F_{u(i,j)},
(\widehat E/\widehat F_{u(i,1)})\widehat F_{u(i,1)Y})\\
&=\text{int}
\left(\widehat F_{u(i,j)},\sum_{0\le i'\le w\text{ with }i'\ne i}
(\widehat E/\widehat F_{u(i,1)})\widehat F_{u(i',1)Y}\right)\\
&=\text{int}(\widehat F_{u(i,j)},\widehat E_Y),
\end{align*}
we get
\begin{align*}
\text{int}(\widehat F_{u(i,j)},J(\widehat E,\widehat G))
=&\text{int}(\widehat F_{u(i,j)},\widehat G)
+\text{int}(\widehat F_{u(i,j)},\widehat E_Y)\\
&-\widehat N_{u(i,j)}+\beta(\widehat F_{u(i,j)},\widehat G)
\end{align*}
where each term is an integer.
Summing the above equation over $1\le i\le w$ and $1\le j\le v(i)$, 
and noting that
$$
\sum_{1\le i\le w}\sum_{1\le j\le v(i)}
\text{int}(\widehat F_{u(i,j)},J(\widehat E,\widehat G))
=\text{int}(\widehat F,J(\widehat E,\widehat G))
$$
and
$$
\sum_{1\le i\le w}\sum_{1\le j\le v(i)}
\text{int}(\widehat F_{u(i,j)},\widehat G)
=\text{int}(\widehat F,\widehat G)
$$
and
$$
\sum_{1\le i\le w}\sum_{1\le j\le v(i)}
\text{int}(\widehat F_{u(i,j)},\widehat E_Y)
=\text{int}(\widehat F,\widehat E_Y)
$$
and
$$
\sum_{1\le i\le w}\sum_{1\le j\le v(i)}\widehat N_{u(i,j)}=N
$$
and
$$
\sum_{1\le i\le w}\sum_{1\le j\le v(i)}
\beta(\widehat F_{u(i,j)},\widehat G)
=\beta(\widehat F,\widehat G),
$$
we get
$$
\text{int}(\widehat F,J(\widehat E,\widehat G))
=\text{int}(\widehat F,\widehat G)
+\text{int}(\widehat F,\widehat E_Y)
-N+\beta(\widehat F,\widehat G)
$$
where each term is an integer.
Now the desired identity follows from the above equation by noting that
clearly
$$
\text{int}(\widehat F,\widehat G)=\text{int}(F,G)\nu
\quad\text{ and }\quad
\text{int}(\widehat F,\widehat E_Y)=\text{int}(F,E_Y)\nu
$$
and
$$
\beta(\widehat F,\widehat G)=\beta(F,G)\nu
$$
and, upon letting $\widehat J(X,Y)$ be obtained by substituting $X^{\nu}$ 
for $X$ in the expression of $J(E,G)$ as a member of $k((X))[Y]$, by
the chain rule for jacobians we have
$$
J(\widehat E,\widehat G)=\widehat J(X,Y)J(X^{\nu},Y)
=\widehat J(X,Y)\nu X^{\nu -1}
$$
and hence
$$
\text{int}(\widehat F,J(\widehat E,\widehat G))
=\text{int}(F,J(E,G))\nu+N(\nu-1).
$$

Now we shall deduce:

\centerline{}

{\bf The beta-derivative identity (2.2).} {\it Let $F\in R$ be $k$-monic of
$Y$-degree $N$, and let $G\in R$ be such that {\rm GCD}$(F,G-c)=1$ for all 
$c\in k$ and, in the notation of Section 1, 
$G_Y\in F_jR$ for $1\le j\le\chi(F)$. Then we have
$$
{\rm int}(F,G_X)={\rm int}(F,G)-N+\beta(F,G)
$$
where each term is an integer.}

\centerline{}

PROOF. Let $E=\text{rad}(F)$ and recall that
$J(E,G)=E_XG_Y-E_YG_X$; since by assumption $G_Y\in F_jR$ for 
$1\le j\le\chi(F)$, we get
$\text{int}(F,J(E,G))=\text{int}(F,E_Y)+\text{int}(F,G_X)$; therefore
the beta-derivative indentity follows by subtracting int$(F,E_Y)$ from
both sides of the beta-jacobian identity.

\centerline{}

Next we shall deduce:

\centerline{}

{\bf The betabar-derivative identity (2.3).} {\it Let $F\in R$ be $k$-monic of
$Y$-degree $N>0$ with {\rm GCD}$(F_Y,F-c)=1$ for all $c\in k$. Then we have
$$
{\rm int}(F_X,F_Y)={\rm int}(F,F_Y)-N+1+\overline\beta(F)
$$
where each term is an integer.}

\centerline{}

PROOF.  This follows by taking $(F_Y,F)$ for $(F,G)$ in the beta-derivative
identity.

\centerline{}

{\bf Remark (2.4).} 
Let the notation be as in Section 1, with $F$ and $G$ in $R$.

Note that for generic $\lambda\in k$, i.e., for all except a finite
number of $\lambda\in k$ we have 
$$
\text{int}(F_j,G-\lambda)\le 0\text{ for all }j.
$$
This motivates calling $G$ {\bf generic} (relative to $F$) to mean that
$$
\text{int}(F_j,G)\le 0\text{ for all }j.
$$
It is clear that
\begin{equation*}
G-\lambda\text{ is generic for almost all }\lambda\in k
\tag{2.4.1}
\end{equation*}
(where almost all means for all except a finite number of) and
\begin{equation*}
G\text{ is generic }\Leftrightarrow\text{int}(F,G)=\text{minint}(F,G)
\tag{2.4.2}
\end{equation*}
and
\begin{equation*}
G\text{ generic }
\Leftrightarrow\text{int}(F,G)=\text{minint}(F,G)
\Leftrightarrow 0\not\in\alpha(F,G)
\tag{2.4.3}
\end{equation*}
and
\begin{equation*}
G\text{ generic }\Rightarrow
\beta(F,G)=\sum_{\lambda\in k}[\text{int}(F,G-\lambda)-\text{int}(F,G)]
\tag{2.4.4}
\end{equation*}
which motivates calling $\beta(F,G)$ the relative excess intersection of $G$.

\centerline{}

\centerline{{\bf Section 3: The conductor-derivative formula}}

\centerline{}

Let $R_0=k[[X,Y]]$.
We say that $F=F(X,Y)\in R$ is {\bf $k$-distinguished} if it is $k$-monic of
$Y$-degree $N$ and, for $0\le i<N$, the coefficient of $Y^i$ in
it has positive $X$-order; note that then $0\ne F\in R_0$ and: 
$N>0\Leftrightarrow F(0,0)=0$. For any $0\ne F=F(X,Y)\in R_0$ with $F(0,0)=0$,
we let $B(F)=R_0/(FR_0)$ and call it the {\bf local ring} of $F$; clearly 
$B(F)$ is a one-dimensional local ring.
Recall that the conductor $C$ of
a ring $B$ is the largest ideal in $B$ which remains an ideal in the
integral closure $B^*$ of $B$ in its total quotient ring;
the length of $C$ is the maximum length $n$ of strictly increasing
chains of ideals $C=C_0<C_1<\dots<C_n\subset B$ in $B$; if there is no
maximum we take the length to be $\infty$; otherwise it is a nonnegative
integer; note that $n=0\Leftrightarrow C=B\Leftrightarrow B=B^*$.
For $0\ne F=F(X,Y)\in R_0$ with $F(0,0)=0$ 
we let $\delta(F)$ denote the length of the conductor of $B(F)$
and call it the {\bf conductor-length} of $F$; 
note that if $F$ is $k$-distinguished with rad$(F)=F$ then $\delta(F)$ is a 
nonnegative integer. 

We shall now prove two lemmas and then we shall prove the 
conductor-derivative formula.

\centerline{}

{\bf Lemma (3.1).} {\it Let $F=F(X,Y)\in R_0$ be $k$-distinguished
of $Y$-degree $N\ge 0$ and let $G=G(X,Y)\in R_0$ be such that $G(0,0)=0$. 
Then $\alpha(F,G)\subset\{0\}$ and $\beta(F,G)=0$.}

\centerline{}

PROOF.
In the notation of Section 1, $F$ distinguished $\Rightarrow z_i(X)\in k[[X]]$
with $z_i(0)=0$ for $1\le i\le N$, and hence 
$\alpha(F,G)\subset\{0\}$ and therefore $\beta(F,G)=0$.

\centerline{}

{\bf Lemma (3.2).} {\it Let $F=F(X,Y)\in R_0$ be $k$-distinguished
of $Y$-degree $N>0$. Then
$\overline\alpha(F)\subset\{0\}$ and $\overline\beta(F)=0$.}

\centerline{}

PROOF.
Follows from (3.1).

\centerline{}

{\bf The conductor-derivative formula (3.3).} {\it Let $F\in R_0$ be
$k$-distinguished of $Y$-degree $N>0$ with {\rm rad}$(F)=F$. Then
\begin{equation*}
{\rm int}(F_X,F_Y)={\rm int}(F,F_Y)-N+1
\tag{3.3.1}
\end{equation*}
and
\begin{equation*}
{\rm int}(F,F_Y)-N+1=2\delta(F)-\chi(F)+1
\tag{3.3.2}
\end{equation*}
and
\begin{equation*}
{\rm int}(F_X,F_Y)=2\delta(F)-\chi(F)+1
\tag{3.3.3}
\end{equation*}
where all the terms in the above three equations are integers.}

\centerline{}

PROOF.
In view of (3.2), the first equation follows from the betabar-derivative 
identity. The second equation follows from Dedekind's Theorem; see pages
65 and 150 of \cite{Ab4}. 
The third equation follows from the first two.

\centerline{}

\centerline{{\bf Section 4: The affine beta-jacobian identity}}

\centerline{}

Let $R_2$ be the affine coordinate ring of the affine plane over $k$, i.e.,
let $R_2$ be the polynomial ring $k[X,Y]$, and note that then 
$$
R_2=k[X,Y]\subset K[[X]][Y]=R_0\cap R\subset K((X))((Y)).
$$
We are particularly interested in the case when $F(X,Y)=f(X^{-1},Y)$
with $f(X,Y)$ in $R_2$; when this is so we shall say that $F$ is the 
{\bf meromorphic associate} of $f$ or $f$ is the {\bf polynomial associate} 
of $F$, and we indicate it by writing $F\sim_m f$ or $f\sim_p F$.
Likewise when $F(X,Y)=f(X^{-1},Y)$ and $G(X,Y)=g(X^{-1},Y)$
we shall indicate it by writing $(F,G)\sim_m(f,g)$ or $(f,g)\sim_p(F,G)$.

Referring for details to Appendix I of Abhyankar's Montreal Lecture Notes
\cite{Ab1} and Lectures 5 to 19 of Abhyankar's Engineering Book \cite{Ab4}, 
let us review some definitions and facts about the
{\bf intersection multiplicity} int$(f,g;\mathcal A)$ of two plane curves
$f$ and $g$ in the {\bf affine plane} $\mathcal A=k^2$ over $k$, i.e., of
two members $f=f(X,Y)$ and $g=g(X,Y)$ of $R_2$.
Let
$$
N=\text{deg}_{(X,Y)}f\quad\text{ and }\quad M=\text{deg}_{(X,Y)}g
$$
where deg$_{(X,Y)}f$ denotes the $(X,Y)$-degree of $f$, i.e., the total
degree of $f$; note that $N$ is a nonnegative integer or $-\infty$ according
as $f\ne 0$ or $f=0$.
Let us write gcd$(f,g)=1$ or gcd$(f,g)\ne 1$ according as the curves $f$ 
and $g$ do not or do have a {\bf common component}, i.e., according
as the polynomials $f$ and $g$ do not or do have a {\bf nonconstant} 
(= not belonging to $k$) common factor in $R_2$. 
If gcd$(f,g)\ne 1$ then we put int$(f,g;\mathcal A)=\infty$.
If gcd$(f,g)=1$ then the curves $f$ and $g$ have a finite 
number of common points $Q_i=(u_i,v_i)$ for $1\le i\le n$ in $\mathcal A$, 
i.e., $f(u_i,v_i)=0=g(u_i,v_i)$ for $1\le i\le n$, and we put
$$
\text{int}(f,g;\mathcal A)=\sum_{1\le i\le n}\text{int}(f,g;Q_i)
$$
where for any $Q=(u,v)\in\mathcal A$ we define the 
{\bf intersection multiplicity} int$(f,g;Q)$ of $f$ and $g$ at $Q$ thus.
Given any $Q=(u,v)\in\mathcal A$, if $f\ne 0$ then by
the Weierstrass Preparation Theorem we can uniquely write
$$
f(X+u,Y+v)=\widehat f_Q(X,Y)f_Q(X,Y)
$$
where
$$
\widehat f_Q(X,Y)\in k[[X,Y]]\text{ with }\widehat f_Q(0,0)\ne 0
$$
and
$$
f_Q(X,Y)=X^a[Y^b+c_1(X)Y^{b-1}+\dots+c_b(X)]
$$
with nonnegative integers $a,b$ and elements $c_1(X),\dots,c_b(X)$ in $k[[X]]$
for which $c_1(0)=\dots=c_b(0)=0$, and we call $f_Q=f_Q(X,Y)\in k[[X,Y]]$ the 
{\bf incarnation} of $f$ at $Q$; if $f=0$ then we put $f_Q=f_Q(X,Y)=0$.
Now for any $Q\in\mathcal A$ we define
$$
\text{int}(f,g;Q)=\text{int}(f_Q,g_Q).
$$
It can be shown that in terms of $k$-vector space dimension we always have
\begin{equation*}
\text{int}(f,g;\mathcal A)=[R_2/(f,g)R_2:k].
\tag{4.1}
\end{equation*}

To extend the above discussion to the {\bf intersection multiplicity}
int$(f,g;\mathcal P)$ of
$f$ and $g$ in the {\bf projective plane} $\mathcal P$ over $k$, recall that
$\mathcal P$ is the disjoint union of
$\mathcal A$ with the line at infinity $\mathcal L_{\infty}$ given by
$$
\mathcal L_{\infty}=\{(\infty,v):v\in k\cup\{\infty\}\}.
$$
If $f\ne 0$ then we define its
{\bf homogenization} $f_h=f_h(X,Y,Z)\in k[X,Y,Z]$ by putting
$$
f_h(X,Y,Z)=Z^N f(X/Z,Y/Z)
$$
and if $f=0$ then we put $f_h=f_h(X,Y,Z)=0$.
For any $Q=(\infty,b)\in \mathcal L_{\infty}$ we define the {\bf incarnation} 
$f_Q=f_Q(X,Y)\in k[[X,Y]]$ of $f$ at $Q$ by putting
$$
f_Q(X,Y)=\begin{cases}
(f_h(1,Y,X))_{(0,b)}&\text{if }b\in k\\
(f_h(X,1,Y))_{(0,0)}&\text{if }b=\infty\\
\end{cases}
$$
and we define int$(f,g;Q)$ by putting
$$
\text{int}(f,g;Q)=\text{int}(f_Q,g_Q).
$$
Now $Q=(u,v)\in\mathcal A$ is a common point of $f$ and $g$, i.e.,
$f(u,v)=0=g(u,v)$, iff (= if and only if) $f_Q(0,0)=0=g_Q(0,0)$; by analogy 
we call $Q\in\mathcal L_{\infty}$ a {\bf common point} of $f$ and $g$
iff $f_Q(0,0)=0=g_Q(0,0)$.
If gcd$(f,g)\ne 1$ then we put int$(f,g;\mathcal P)=\infty$.
If gcd$(f,g)=1$ then the curves $f$ and $g$ have a finite 
number of common points $Q_i=(u_i,v_i)$ for $1\le i\le m$ in $\mathcal P$, 
i.e., $f_{Q_i}(0,0)=0=g_{Q_i}(0,0)$ for $1\le i\le m$, and we put
$$
\text{int}(f,g;\mathcal P)=\sum_{1\le i\le m}\text{int}(f,g;Q_i).
$$
Note that by {\bf Bezout's Theorem}
\begin{equation*}
\text{int}(f,g;\mathcal P)=\begin{cases}
NM&\text{if gcd$(f,g)=1$ and $f\ne 0\ne g$}\\
0&\text{if gcd$(f,g)=1$ and either $f=0\ne g$ or $f\ne 0=g$}\\
\infty&\text{if gcd$(f,g)\ne 1$.}\\
\end{cases}
\tag{4.2}
\end{equation*}

For simplicity of notation let us put
$$
(\infty)=(\infty,0)\in\mathcal L_{\infty}.
$$
We say that $f$ is {\bf $Y$-monic} to mean that 
$$
f\ne 0\text{ and deg}_{(X,Y)}[f(X,Y)-f_0Y^N]<N\text{ for some }0\ne f_0\in k
$$
and we note that 
$$
f\text{ is $Y$-monic }\Rightarrow\text{ the $(X,Y)$-degree of $f$ coincides
with its $Y$-degree}
$$
and
$$
f\text{ is $Y$-monic }\Leftrightarrow 
\{Q\in\mathcal L_{\infty}:f_Q(0,0)=0\}\subset\{(\infty)\}.
$$
It can be shown that
\begin{equation*}
\begin{cases}
\text{if $(F,G)\sim_m(f,g)$ where $f$ is $Y$-monic and gcd$(f,g)=1$,}\\
\text{then int}(F,G)=-\text{int}(f,g;\mathcal A).
\end{cases}
\tag{4.3}
\end{equation*}

If $f\ne 0$ and $0\ne f_0\in k[X]$ is the coefficient of the highest power 
of $Y$ in $f$, then by the {\bf top coefficient} of $f$ we mean the nonzero
element of $k$ which is the coefficient of the highest power of $X$ in $f_0$;
if $f=0$ then we take $0$ to be the top coefficient of $f$. If $f\in k$ 
then we define the {\bf radical} of $f$ in $\mathcal A$ by putting
rad$(f;\mathcal A)=f$, and if $f\not\in k$ then we factor $f$ as a product of
positive powers of mutually nonassociate irreducible members of $R_2$ by
writing $f=f_1^{e(1)}\dots f_r^{e(r)}$ and we put 
rad$(f;\mathcal A)=cf_1\dots f_r$ where $0\ne c\in k$ is so chosen that
$f$ and rad$(f;\mathcal A)$ have the same top coefficient.

In analogy with the alpha and  beta invariants, 
for any $f,g$ in $R_2$ we define 
the {\bf affine alpha invariant} $\alpha(f,g;\mathcal A)$ and the 
{\bf affine beta invariant} $\beta(f,g;\mathcal A)$ of $f$ relative to $g$,
and in analogy with the alphabar and beta-bar invariants,
for any $f\in R_2$ we define
the {\bf affine alphabar invariant} $\overline\alpha(f;\mathcal A)$ and
the {\bf affine betabar invariant} $\overline\beta(f;\mathcal A)$ of $f$
thus.

First, for any $f,g$ in $R_2$, 
we define the {\bf maximal-intersection} of $f$ with $g$ by putting
$$
\text{maxint}(f,g;\mathcal A)
=\text{max}_{\mu\in k}\text{int}(f,g-\mu;\mathcal A)
$$
(where max stands for lub = least upper bound) and we note that by Bezout
$$
\begin{cases}
\text{maxint}(f,g;\mathcal A)=\text{a nonnegative integer or }\infty,\\
\text{and: maxint}(f,g;\mathcal A)=\infty
\Leftrightarrow\text{gcd}(f,g-c;\mathcal A)\ne 1\text{ for some }c\in k
\end{cases}
$$
and we put
$$
\alpha(f,g;\mathcal A)
=\{\lambda\in k:\text{int}(f,g-\lambda;\mathcal A)
<\text{maxint}(f,g;\mathcal A)\}
$$
and
$$
\beta(f,g;\mathcal A)=\sum_{0\ne\lambda\in\alpha(f,g;\mathcal A)}
[\text{maxint}(f,g;\mathcal A)-\text{int}(f,g;\mathcal A)]
$$
and we note that then
$$
\beta(f,g;\mathcal A)=\text{a nonnegative integer or }\infty.
$$
By a {\bf relative affine irregular value} of $g$ (on $f$) we mean an 
element of $\alpha(f,g;\mathcal A)$. We call $\alpha(f,g;\mathcal A)$ the 
{\bf relative affine irregular value set} of $g$ (on $f$), and we call 
$\beta(f,g;\mathcal A)$ the {\bf relative deficit intersection} of 
$g$ (on $f$).

Next, for any $f$ in $R_2$, we put
$$
\overline\alpha(f;\mathcal A)=\alpha(f_Y,f;\mathcal A)
$$ 
and 
$$
\overline\beta(f;\mathcal A)=\beta(f_Y,f;\mathcal A).
$$
By an {\bf affine irregular value} of $f$ we mean an element of
$\overline\alpha(f;\mathcal A)$. We call $\overline\alpha(f;\mathcal A)$ the 
{\bf affine irregular value set} of $f$, and we call 
$\overline\beta(f;\mathcal A)$ the {\bf deficit intersection} of $f$. 

As a consequence of (4.3) we see that:
\begin{equation*}
\begin{cases}
\text{if $(F,G)\sim_m(f,g)$ where $f$ is $Y$-monic}\\
\text{and gcd$(f,g-c)=1$ for all $c\in k$,}\\
\text{then $\alpha(f,g;\mathcal A)=\alpha(F,G)=$ a finite set}\\
\text{and $\beta(f,g;\mathcal A)=\beta(F,G)=$ a nonnegative integer.}\\
\end{cases}
\tag{4.4}
\end{equation*}
and
\begin{equation*}
\begin{cases}
\text{if $F\sim_m f$ where $f$ is $Y$-monic of $Y$-degree $N>0$}\\
\text{and gcd$(f_Y,f-c)=1$ for all $c\in k$,}\\
\text{then $\overline\alpha(f;\mathcal A)
=\overline\alpha(F)=$ a finite set}\\
\text{and $\overline\beta(f;\mathcal A)=\overline\beta(F)=$ 
a nonnegative integer.}\\
\end{cases}
\tag{4.5}
\end{equation*}

Now as consequences of (2.1), (2.2), (2.3) we shall deduce
(4.6), (4.7), (4.8) respectively.

\centerline{}

{\bf The affine beta-jacobian identity (4.6).} 
{\it Let $f\in R_2$ be $Y$-monic of $Y$-degree $N$, and let $g\in R_2$ be 
such that {\rm gcd}$(f,g-c)=1$ for all $c\in k$. 
Then for $e={\rm rad}(f;\mathcal A)$ we have
$$
{\rm int}(f,g;\mathcal A)+{\rm int}(f,e_Y;\mathcal A)
={\rm int}(f,J(e,g);\mathcal A)+\beta(f,g;\mathcal A)+N
$$
where each term is a nonnegative integer.}

\centerline{}

PROOF. Let $(F,G)\sim_m(f,g)$; then $F\in R$ is $k$-monic of $Y$-degree
$N$, and $G\in R$ is such that GCD$(F,G-c)=1$ for all $c\in k$.
Let $E\sim_m e$; then $E=\text{rad}(F)$ and $E_Y\sim_m e_Y$.
Let $J=J(E,G)$; also let $j'=J(e,g)$ and $J'\sim_m j'$;
then $J'$ is the jacobian of $E$ and $G$ relative to $X^{-1}$ and $Y$
and hence by the chain rule for jacobians we see that $J$ equals $J'$
times the jacobian of $X^{-1}$ and $Y$ relative to $X$ and $Y$ which is
clearly equal to $-X^{-2}$. Also obviously
int$(F,-X^{-2})=-2N$ and 
int$(F,J)=\text{int}(F,-X^{-2}J')
=\text{int}(F,-X^{-2})+\text{int}(F,J')$. Therefore
$\text{int}(F,J)=-2N+\text{int}(F,J')$, 
and substituting this in (2.1) we get
$$
-2N+\text{int}(F,J')=\text{int}(F,G)+\text{int}(F,E_Y)-N+\beta(F,G)
$$
and by sending some terms from one side to the other we obtain
$$
-\text{int}(F,G)-\text{int}(F,E_Y)=-\text{int}(F,J')+\beta(F,G)+N
$$
and in view of (4.3) and (4.4) this yields the desired result.

\centerline{}

{\bf The affine beta-derivative identity (4.7)} {\it Let $f\in R_2$ be 
$k$-monic of $Y$-degree $N$, and let $g\in R_2$ be such that 
{\rm gcd}$(f,g-c)=1$ for all $c\in k$ and, in the notation of Section 1, 
for $(F,G)\sim_m(f,g)$ we have
$G_Y\in F_jR$ for $1\le j\le\chi(F)$. Then 
$$
{\rm int}(f,g;\mathcal A)={\rm int}(f,g_X;\mathcal A)+\beta(f,g;\mathcal A)+N
$$
where each term is a nonnegative integer.}

\centerline{}

PROOF. Now $F\in R$ is $k$-monic of $Y$-degree
$N$, and $G\in R$ is such that GCD$(F,G-c)=1$ for all $c\in k$ and, 
in the notation of Section 1, $G_Y\in F_jR$ for $1\le j\le\chi(F)$.
Let $g'=g_X$ and $G'\sim_m g'$; then $G'$ is the derivative of $G$ 
with respect to $X^{-1}$
and hence by the chain rule for derivatives we see that $G_X$ equals $G'$
times the  derivative of $X^{-1}$ with respect to $X$ which is
clearly equal to $-X^{-2}$. Also obviously
int$(F,-X^{-2})=-2N$ and 
int$(F,G_X)=\text{int}(F,-X^{-2}G')
=\text{int}(F,-X^{-2})+\text{int}(F,G')$. Therefore
$\text{int}(F,G_X)=-2N+\text{int}(F,G')$, 
and substituting this in (2.2) we get
$$
-2N+\text{int}(F,G')=\text{int}(F,G)-N+\beta(F,G)
$$
and by sending some terms from one side to the other we obtain
$$
-\text{int}(F,G)=-\text{int}(F,G')+\beta(F,G)+N
$$
and in view of (4.3) and (4.4) this yields the desired result.

{\bf The affine betabar-derivative identity (4.8).} 
{\it Let $f\in R_2$ be $Y$-monic of $Y$-degree $N>0$ with 
{\rm gcd}$(f_Y,f-c)=1$ for all $c\in k$. Then we have
$$
{\rm int}(f,f_Y;\mathcal A)
={\rm int}(f_X,f_Y;\mathcal A)+\overline\beta(f;\mathcal A)
+(N-1)
$$
where each term is a nonnegative integer.}

PROOF. Let $F\sim_m f$; then $F\in R$ is $k$-monic of $Y$-degree
$N>0$ with GCD$(F,G-c)=1$ for all $c\in k$; also clearly
$F_Y\sim_m f_Y$.
Let $f'=f_X$ and $F'\sim_m f'$; then $F'$ is the derivative of $F$ 
with respect to $X^{-1}$
and hence by the chain rule for derivatives we see that $F_Y$ equals $F'$
times the  derivative of $X^{-1}$ with respect to $X$ which is
clearly equal to $-X^{-2}$. Also obviously
int$(-X^{-2},F_Y)=-2(N-1)$ and 
int$(F_X,F_Y)=\text{int}(-X^{-2}F',F_Y)
=\text{int}(-X^{-2},F_Y)+\text{int}(F',F_Y)$. Therefore
$\text{int}(F_X,F_Y)=-2(N-1)+\text{int}(F',F_Y)$, 
and substituting this in (2.3) we get
$$
-2(N-1)+\text{int}(F',F_Y)=\text{int}(F,F_Y)-N+1+\overline\beta(F)
$$
and by sending some terms from one side to the other we obtain
$$
-\text{int}(F,F_Y)=-\text{int}(F',F_Y)+\overline\beta(F)+(N-1)
$$
and in view of (4.3) and (4.5) this yields the desired result.

\centerline{}

As a consequence of (4.8) we shall now prove:

\centerline{}

{\bf The projective betabar-derivative identity (4.9).} 
{\it Let $f\in R_2$ be $Y$-monic of $Y$-degree $N>0$ with 
{\rm gcd}$(f_Y,f-c)=1$ 
for all $c\in k$. Then we have
$$
{\rm int}(f,f_Y;\mathcal P)
={\rm int}(f_X,f_Y;\mathcal P)+\overline\beta(f;\mathcal A)+2(N-1)
$$
where each term is a nonnegative integer.}

PROOF. By (3.2) we have $\overline\beta(f_{(\infty)})=0$ and hence by
the betabar-derivative identity (2.3) we get
$$
\text{int}(f_{(\infty)},f_{(\infty)Y})
=\text{int}(f_{(\infty)X},f_{(\infty)Y})+(N-1)
$$
and adding this to (4.8) we obtain the desired result by noting that
obviously $f_Y$ is $Y$-monic and hence
\begin{equation*}
\text{int}(f,f_Y;\mathcal A)+\text{int}(f_{(\infty)},f_{(\infty)Y})
=\text{int}(f,f_Y;\mathcal P)
\tag{4.9.1}
\end{equation*}
and
\begin{equation*}
\text{int}(f_X,f_Y;\mathcal A)+\text{int}(f_{(\infty)X},f_{(\infty)Y})
=\text{int}(f_X,f_Y;\mathcal P).
\tag{4.9.2}
\end{equation*}

\centerline{}

Finally as a consequence of (4.9) we shall prove:

\centerline{}

{\bf The betabar-conductor identity (4.10).} 
{\it Let $f\in R_2$ be $Y$-monic of $Y$-degree $N>0$ with 
{\rm gcd}$(f_Y,f-c)=1$ 
for all $c\in k$. Then we have
$$
(N-1)(N-2)+[\chi(f_{(\infty)})-1]
={\rm int}(f_X,f_Y;\mathcal A)
+2\delta(f_{(\infty)})+\overline\beta(f;\mathcal A)
$$
where each term is a nonnegative integer.}

\centerline{}

NOTE. We regard $(N-1)(N-2)$ as a term. Also we regard a square bracketed
expression as a term. For instance $[\chi(f_{(\infty)})-1]$ is a term.

\centerline{}

PROOF. By Bezout we have int$(f,f_Y;\mathcal P)=N(N-1)$, and hence by (4.9)
and (4.9.2) we get
$$
(N-1)(N-2)
=\text{int}(f_X,f_Y;\mathcal A)+\text{int}(f_{(\infty)X},f_{(\infty)Y})
+\overline\beta(f;\mathcal A)
$$
and from this our assertion follows by using (3.3.3) applied to
$F=f_{(\infty)}$.

\centerline{}

\centerline{{\bf Section 5: No Irregular Value}}

\centerline{}

For nonconstant $f\in R_2$ we let $B(f;\mathcal A)=R_2/(fR_2)$ and call it the 
{\bf affine coordinate ring} of $f$. Clearly $(fR_2)\cap k=\{0\}$ and hence
we may identify $k$ with a subfield of $B(f;\mathcal A)$. We define
$\delta(f;\mathcal A)$ and $\delta(f;\mathcal P)$ by putting
$$
\delta(f;\mathcal A)=\sum_{Q\in\mathcal A}\delta(f_Q)
$$
and
$$
\delta(f;\mathcal P)=\sum_{Q\in\mathcal P}\delta(f_Q)
$$
and we note that if rad$(f;\mathcal A)=f$ then these are nonnegative integers;
otherwise they are $\infty$. 
In case $f$ is 
irreducible (in $R_2$), the {\bf genus} of the quotient field of 
$B(f;\mathcal A)$ over $k$ is a nonnegative integer which we denote by 
$\gamma(f)$.  Again if $f$ is irreducible then by the 
{\bf number of places of $f$ at infinity}, 
denoted by $\chi_{\infty}(f)$, we mean the number of DVRs
(= discrete valuations rings = one dimensional regular local rings) 
which do not contain $B(f;\mathcal A)$ but whose
quotient field coincides with the quotient field of $B(f;\mathcal A)$.
In the general case, by factoring $f$ as a product of irreducible
factors, i.e., writing $f=f_1\dots f_r$ where $f_1,\dots,f_r$ 
are irreducible members of $R_2$, we call
$\chi_{\infty}(f_1)+\dots+\chi_{\infty}(f_r)$ the 
{\bf number of places of $f$ at infinity}, and denote it by
$\chi_{\infty}(f)$; note that this is always a positive integer;
also note that
$$
\begin{cases}
\text{if $f$ is $Y$-monic then $\chi_{\infty}(f)=\chi(f_{(\infty)})$ and}\\
\text{upon letting $F\sim_m f$ we have $\chi(F)=\chi_{\infty}(f)$.}
\end{cases}
$$
Let us call $f$ a {\bf uniline} if $f$ is irreducible with
$\gamma(f)=0=\delta(f;\mathcal A)$ and $\chi_{\infty}(f)=1$, and
let us call $f$ a {\bf unihyperbola} if $f$ is irreducible with
$\gamma(f)=0=\delta(f;\mathcal A)$ and $\chi_{\infty}(f)=2$. Also let us call
$f$ a {\bf multihyperbola} if in the above notation we have
int$(f_i,f_j;\mathcal A)=0$ for $1\le i<j\le r$ and $f_i$ is a unihyperbola
for $1\le i\le r$, and let us call
$f$ a {\bf multihyperbolic line} if in the above notation we have
int$(f_i,f_j;\mathcal A)=0$ for $1\le i<j\le r$, $f_j$ is a uniline
for some $j\in\{1,\dots,r\}$, and $f_i$ is a unihyperbola for all
$i\in\{1,\dots,j-1,j+1,\dots,r\}$.
Clearly
\begin{align*}
f\text{ is a uniline }&\Leftrightarrow B(f;\mathcal A)\approx k[X]\\
&\Leftrightarrow f\text{ is an irreducible multihyperbolic line}
\end{align*}
and
$$
f\text{ is a unihyperbola }\Leftrightarrow B(f;\mathcal A)\approx k[X,X^{-1}]
$$
where $\approx$ denotes $k$-isomorphism of rings.
By the {\bf epimorphism theorem} (see \cite{Ab2}) we know that
\begin{equation*}
\begin{cases}
f\text{ is a uniline}\\
\Leftrightarrow\text{ for some $k$-automorphism $\sigma$ of $R_2$
we have $\sigma(f)=X$}\\
\end{cases}
\tag{5.1}
\end{equation*}
and
\begin{equation*}
\begin{cases}
f\text{ is a multihyperbolic line which is not a uniline}\\
\Leftrightarrow\text{ for some $k$-automorphism $\sigma$ of $R_2$}\\
\qquad\text{we have $\sigma(f)=X(1+Xh(X,Y))$}\\
\qquad\text{where $0\ne h(X,Y)\in R_2$ is such that}\\ 
\qquad\text{$1+Xh(X,Y)$ is a multihyperbola.}
\end{cases}
\tag{5.2}
\end{equation*}
The well-known {\bf genus formula} tells us that
\begin{equation*}
\begin{cases}
\text{if $f\in R_2$ is irreducible of $(X,Y)$-degree $N>0$ then}\\
2\gamma(f)+2\delta(f;\mathcal P)=(N-1)(N-2).
\end{cases}
\tag{5.3}
\end{equation*}
Likewise, continuing the assumption of nonconstant $f\in R_2$,
a well-known consequence of {\bf Bertini's Theorem} tells us that
\begin{equation*}
\begin{cases}
\text{rad$(f-c;\mathcal A)=f-c$ for all $c\in k$}\\
\text{$\Rightarrow f-\lambda$ is irreducible for almost all $\lambda\in k$}
\end{cases}
\tag{5.4}
\end{equation*}
and as supplements to this we note the easy to prove facts which say that
\begin{equation*}
\begin{cases}
\text{rad$(f-c;\mathcal A)=f-c$ for all $c\in k$}\\
\Leftrightarrow\text{int}(f_X,f_Y;\mathcal A)<\infty\\
\end{cases}
\tag{5.5}
\end{equation*}
and
\begin{equation*}
\begin{cases}
\text{int}(f_X,f_Y;\mathcal A)=0\Rightarrow\text{int}(f_X,f_Y;\mathcal A)<\infty,\\
\text{and}\\
\text{int}(f_X,f_Y;\mathcal A)=0
\Leftrightarrow\delta(f-c;\mathcal A)=0\text{ for all }c\in k\\
\end{cases}
\tag{5.6}
\end{equation*}
and
\begin{equation*}
\begin{cases}
\text{if $f$ is  $Y$-monic then:}\\
\text{rad$(f-c;\mathcal A)=f-c$ for all $c\in k$}\\
\text{$\Leftrightarrow$ gcd$(f_Y,f-c;\mathcal A)=1$ for all $c\in k$.}
\end{cases}
\tag{5.7}
\end{equation*}
We shall now deduce some consequences of (5.1) to (5.7).

\centerline{}

{\bf The betabar-genus identity (5.8).} 
{\it Let $f\in R_2$ be irreducible $Y$-monic of $Y$-degree $N>0$ with 
{\rm gcd}$(f_Y,f-c)=1$ for all $c\in k$. Then
\begin{equation*}
2\gamma(f)+2\delta(f;\mathcal A)+[\chi(f_{(\infty)})-1]
={\rm int}(f_X,f_Y;\mathcal A)+\overline\beta(f;\mathcal A)
\tag{5.8.1}
\end{equation*}
where each term is a nonnegative integer.
Moreover,
\begin{equation*}
\begin{cases}
\text{if {\rm int}$(f_X,f_Y;\mathcal A)=0$}\\
\text{then }2\gamma(f)+[\chi(f_{(\infty)})-1]=\overline\beta(f;\mathcal A)
\end{cases}
\tag{5.8.2}
\end{equation*}
where each term is a nonnegative integer. 
Finally,
\begin{equation*}
\begin{cases}
\text{if {\rm int}$(f_X,f_Y;\mathcal A)=0=\overline\beta(f;\mathcal A)$}\\
\text{then $f$ is a uniline.}
\end{cases}
\tag{5.8.3}
\end{equation*}}

\centerline{}

PROOF. (5.8.1) follows from (4.10) and (5.3) by noting that
$\delta(f;\mathcal P)=\delta(f;\mathcal A)+\delta(f_{(\infty)})$.
Moreover, (5.8.2) follows from (5.8.1) by noting that
int$(f_X,f_Y;\mathcal A)=0\Rightarrow\delta(f;\mathcal A)=0$.
Since $\gamma(f)$ and $[\chi(f_{(\infty)})-1]$ are nonnegative integers,
it follows that if $2\gamma(f)+[\chi(f_{(\infty)})-1]=0$ then we must have
$\gamma(f)=0=[\chi(f_{(\infty)})-1]$; therefore, since
int$(f_X,f_Y;\mathcal A)=0\Rightarrow\delta(f;\mathcal A)=0$, 
(5.8.3) follows from (5.8.2).

\centerline{}

{\bf No Irregular Value Theorem (5.9).}
{\it Let $f\in R_2$ be $Y$-monic of $Y$-degree $N>0$. Then
\begin{align*}
{\rm int}(f_X,f_Y;\mathcal A)=0=|\overline\alpha(f;\mathcal A)|
&\Leftrightarrow f\text{ is a uniline.} \\
&\Leftrightarrow f-c\text{ is a uniline for all }c\in k.\\
\end{align*}}

\centerline{}

PROOF. We shall give a circular proof by showing that
LHS $\Rightarrow$ MHS $\Rightarrow$ RHS $\Rightarrow$ LHS where
MHS = Middle Hand Side = the condition ``$f$ is a uniline.''
Assuming int$(f_X,f_Y;\mathcal A)=0$, 
by (5.4) to (5.7) there exists $\lambda\in k$ such that for
$f'=f-\lambda$ we have that $f'\in R_2$ is irreducible $Y$-monic of
$Y$-degree $N>0$ such that int$(f'_X,f'_Y)=0$ and
gcd$(f'_Y,f'-c)=1$ for all $c\in k$; also assuming 
$|\overline\alpha(f;\mathcal A)|=0$, we see that
$\overline\beta(f';\mathcal A)=0$, and hence $f'$ is a uniline by (5.8.3), and
therefore $f$ is a uniline by (5.1).
By (5.1) we know that if $f$ is a uniline then $f-c$ is a uniline for
all $c\in k$. Finally, upon letting $f_c=f-c$, assume that $f_c$ is a 
uniline for all $c\in k$. Then, for any $c\in k$, clearly
$f_c\in R_2$ is irreducible $Y$-monic of $Y$-degree $N>0$ with
$\gamma(f_c)=0=[\chi(f_{c(\infty)})-1]$, and by (5.5) to (5.7) 
we see that int$(f_{cX},f_{cY})=0$ and gcd$(f_{cY},f_c-c')=1$ for all
$c'\in k$, and hence by (5.8.2) we get
$\overline\beta(f_c;\mathcal A)=0$. Therefore 
$\text{int}(f_X,f_Y;\mathcal A)=0=|\overline\alpha(f;\mathcal A)|$.

\centerline{}

{\bf Remark (5.10).} 
As hinted in the above proof, from the definitions of the affine invariants
$\overline\alpha,\overline\beta,\alpha,\beta$ it immediately follows that
for any $f$ in $R_2$ we have
$$
\overline\alpha(f;\mathcal A)\subset\{0\}
\Leftrightarrow\overline\beta(f;\mathcal A)=0
$$
and
$$
\overline\alpha(f;\mathcal A)=\emptyset
\Leftrightarrow\overline\beta(f-c;\mathcal A)=0\text{ for all }c\in k
$$
and
for any $f,g$ in $R_2$ we have
$$
\alpha(f,g;\mathcal A)\subset\{0\}
\Leftrightarrow\beta(f,g;\mathcal A)=0
$$
and
$$
\alpha(f,g;\mathcal A)=\emptyset
\Leftrightarrow\beta(f,g-c;\mathcal A)=0\text{ for all }c\in k.
$$

\centerline{}

\centerline{\bf Section 6: One Irregular Value}

\centerline{}

Here are some consequences of (5.1) to (5.3) and (5.5) to (5.7) 
which do not use (5.4).

\centerline{}

{\bf Product Identity (6.1).}
{\it Let $f\in R_2$ be $Y$-monic of $Y$-degree $N>0$ with 
{\rm gcd}$(f_Y,f)=1$.
Let $f=\prod_{1\le i\le r}f_i$ be a factorization of $f$ where $f_i\in R_2$
is $Y$-monic of $Y$-degree $N_i>0$ for $1\le i\le r$.
Then 
\begin{align*}
{\rm int}(f,f_Y;\mathcal A)-N
=&2\sum_{1\le i<j\le r}{\rm int}(f_i,f_j;\mathcal A))\\
&+\sum_{1\le i\le r}
[(N_i-1)(N_i-2)-2\delta(f_{i(\infty)})+\chi(f_{i(\infty)})-2]\\
\tag{6.1.1}
\end{align*}
where each term is an integer.
Moreover, 
\begin{equation*}
\begin{cases}
\text{if {\rm int}$(f_X,f_Y;\mathcal A)=0$ and $f_i$ is irreducible
for $1\le i\le r$}\\
\text{then }{\rm int}(f,f_Y;\mathcal A)-N
=\sum_{1\le i\le r}[2\gamma(f_i)+\chi(f_{i(\infty)})-2]\\
\text{and $\delta(f_i;\mathcal A)=0$ for $1\le i\le r$}
\end{cases}
\tag{6.1.2}
\end{equation*}
where each term is an integer.}

\centerline{}

PROOF. By the product rule for derivatives we get
$$
{\rm int}(f,f_Y;\mathcal A)
=2\sum_{1\le i<j\le r}{\rm int}(f_i,f_j;\mathcal A)
+\sum_{1\le i\le r}{\rm int}(f_i,f_{iY};\mathcal A)
$$
and for $1\le i\le r$ we have
\begin{alignat*}{2}
&\text{int}(f_i,f_{iY};\mathcal A)&&\\
&=N_i(N_i-1)-\text{int}(f_{i(\infty)},f_{i(\infty)Y})
&&\text{ by Bezout}\\
&=N_i(N_i-1)-2\delta(f_{i(\infty)})+\chi(f_{i(\infty)})-N_i
&&\text{ by (3.3.2)}\\
&=(N_i-1)(N_i-2)-2\delta(f_{i(\infty)})+\chi(f_{i(\infty)})+(N_i-2)
&&\text{ by simplifying}\\
\notag
\end{alignat*}
and by substituting this value of int$(f_i,f_{iY};\mathcal A)$ in the RHS of 
the above equation we get (6.1.1) by noting that $N=\sum_{1\le i\le r}N_i$.
Moreover, if int$(f_X,f_Y)=0$ then for $1\le i<j\le r$ we have
int$(f_i,f_j;\mathcal A))=0$, and for $1\le i\le r$ we have
$\delta(f_i;\mathcal A)=0$ and $\delta(f_{i(\infty)})=\delta(f_i;\mathcal P)$, 
and hence (6.1.2) follows from (5.3) and (6.1.1).

\centerline{}

{\bf Product Lemma (6.2).}
{\it Let $f\in R_2$ be $Y$-monic of $Y$-degree $N>0$ with 
{\rm int}$(f,f_Y;\mathcal A)=N-1$ and {\rm int}$(f_X,f_Y;\mathcal A)=0$.
Then $f$ is a multihyperbolic line.}

\centerline{}

PROOF. By (6.1.2) we have
$-1=\sum_{1\le i\le r}[2\gamma(f_i)+\chi(f_{i(\infty)})-2]$ and 
$\delta(f_i;\mathcal A)=0$ for $1\le i\le r$. Therefore, since
$\gamma(f_i)$ and $\chi(f_{i(\infty)})-1$ are nonnegative integers for
$1\le i\le r$, just by numerical considerations, we conclude that there is 
a unique $j\in\{1,\dots,r\}$ such that $\gamma(f_j)=0=\chi(f_{i(\infty)})-1$ 
and $\gamma(f_j)=0=\chi(f_{i(\infty)})-2$ for all
$i\in\{1,\dots,j-1,j+1,\dots,r\}$. Thus $f$ is a multihyperbolic line.

\centerline{}

{\bf Product Theorem (6.3).}
{\it Let $f\in R_2$ be $Y$-monic of $Y$-degree $N>0$ 
such that {\rm int}$(f_X,f_Y;\mathcal A)=0=\overline\beta(f;\mathcal A)$
and {\rm gcd}$(f_Y,f-c;\mathcal A)=1$ for all $c\in k$.
Then {\rm int}$(f,f_Y;\mathcal A)=N-1$ and $f$ is a multihyperbolic line.}

\centerline{}

PROOF. Follows from (4.8) and (6.2).

\centerline{}

{\bf One Irregular Value Theorem (6.4).}
{\it Let $f\in R_2$ be $Y$-monic of $Y$-degree $N>0$ with
{\rm int}$(f_X,f_Y;\mathcal A)$ $=0=|\overline\alpha(f;\mathcal A)|-1$, and 
let $f'=f-\lambda$ with $\{\lambda\}=\overline\alpha(f;\mathcal A)$.
Then $f'$ is a multihyperbolic line which is not a uniline.}

\centerline{}

PROOF. Since int$(f_X,f_Y;\mathcal A)=0$, by (5.5) to (5.7) we see that
$f'\in R_2$ is $Y$-monic of $Y$-degree $N>0$ 
such that int$(f'_X,f'_Y;\mathcal A)=0$
and gcd$(f'_Y,f'-c;\mathcal A)=1$ for all $c\in k$.
Since $|\overline\alpha(f;\mathcal A)|-1=0$ and 
$f'=f-\lambda$ with $\{\lambda\}=\overline\alpha(f;\mathcal A)$, we
also see that $\overline\beta(f';\mathcal A)=0$. Therefore by (6.3) we 
conclude that $f'$ is a multihyperbolic line, and by (5.9) we see that it 
is not a uniline.

\centerline{}

{\bf Remark (6.5).} In the proof of (6.1) we could have tried to use (4.8)
by assuming the condition gcd$(f_Y,f-c)=1$ for all $c\in k$.
tried to use (4.8) in its proof.
This did not work because, as shown by the following example, this 
condition is not inherited by the factors of $f$. Namely, by 
taking $f=f_1f_2$ with $f_1=X^2Y+1$ and $f_2=Y$, we see that
gcd$(f_Y,f-c)=1$ for all $c\in k$ does not imply
gcd$(f_{1Y},f_1-c)=1$ for all $c\in k$.

\centerline{}

\centerline{{\bf Section 7: Two Conjectures}}

\centerline{}

Let
$$
k^{\times}=k\setminus\{0\}
$$
and, as in Abhyankar's previous lectures, let $\nz$ denote
an unspecified element of $k^{\times}$.
Now here are {\bf two meromorphic jacobian conjectures} for $F,G$ in $R$
of $Y$-degrees $N$ and $M$ respectively.

\centerline{}

CONJECTURE I: 
$J(F,G)=\nz X^{-2}\Rightarrow\beta(F,G)=0\not\in\overline\alpha(F)$.

\centerline{}

CONJECTURE II: $J(F,G)=\nz X^{-2}\Rightarrow$ either $M|N$ or $N|M$.

\centerline{}

{\bf Remark (7.1).} Note that Conjecture II 
is not true if we allow the coefficients 
of $F$ and $G$ to be fractional meromorphic series. For example, if
$F=Y^3+(3/2)X^{-1/2}Y\text{ and }G=Y^2+X^{-1/2}$
then 
$J(F,G)=(-3/4)X^{-2}$ but $M=2$ and $N=3$.

\centerline{}

{\bf Remark (7.2).} In Remark (8.9) of the next section we shall show
that both these meromorphic conjectures imply the algebraic jacobian
conjecture which predicts that if $f,g$ in $R_2$ are such that
$J(f,g)=\nz$ then $k[f,g]=R_2$.

\centerline{}

\centerline{\bf Section 8: Thoughts on Conjecture I}

\centerline{}

Before turning to Conjecture I, let us observe that for any $f,g$ in $R_2$
we have
\begin{equation*}
J(f,g)=\nz\text{ and $f$ is a uniline }\Rightarrow k[f,g]=R_2
\tag{8.1}
\end{equation*}
and
\begin{equation*}
J(f,g)=\nz\text{ and $f$ is a multihyperbolic line }
\Rightarrow f\text{ is a uniline}
\tag{8.2}
\end{equation*}
and
\begin{equation*}
k[f,g]=R_2\Rightarrow|\alpha(f,g;\mathcal A)|=0.
\tag{8.3}
\end{equation*}
Out of this (8.1) and (8.2) follow from the parallelness of the
Newton polygons of $f$ anf $g$ proved in \cite{Ab3, Ab5}. The third
assertion (8.3) follows by noting that if $k[f,g]=R_2$ then for every
$c\in k$ we clearly have $k[f,g-c]=R_2$ and hence int$(f,g-c;\mathcal A)=1$,
and therefore $|\alpha(f,g;\mathcal A)|=0$.

\centerline{}

Now let us prove the following two polynomial analogue of (5.8.3).

\centerline{}

{\bf No Deficit Intersection Theorem (8.4).} {\it Let $f\in R_2$ be 
$Y$-monic of $Y$-degree $N>0$, and let $g\in R_2$ be such that 
$J(f,g)=\nz$ and $\beta(f,g;\mathcal A)=0\not\in\overline\alpha(f;\mathcal A)$.
Then $f$ is a uniline, and hence in particular
$|\overline\alpha(f;\mathcal A)|=0=\overline\beta(f;\mathcal A)$ 
and $k[f,g]=R_2$.}

\centerline{}

PROOF. Since $J(f,g)=\nz$, we get 
$
\text{rad}(f;\mathcal A)=f\text{ and int}(f,J(f,g);\mathcal A)=0
$ 
with 
$\text{gcd}(f,g-c;\mathcal A)=1\text{ for all }c\in k.
$ 
Therefore, since $\beta(f,g;\mathcal A)=0$, by (4.6) we obtain
$$
\text{int}(f,g;\mathcal A)+\text{int}(f,f_Y;\mathcal A)=N.
$$
Since $J(f,g)=\nz$, we also get
int$(f_X,f_Y;\mathcal A)=0$ with gcd$(f_Y,f-c;\mathcal A)=1$ for all 
$c\in k$. Therefore by (4.8) we obtain
$$
\text{int}(f,f_Y;\mathcal A)=\overline\beta(f;\mathcal A)+(N-1).
$$
The above two displays tell us that
$$
\text{int}(f,g;\mathcal A)+\overline\beta(f;\mathcal A)=1.
$$
Since int$(f,g;\mathcal A)$ and $\overline\beta(f;\mathcal A)$ are 
nonnegative integers, we must have $\overline\beta(f;\mathcal A)=0$ or $1$.
For a moment suppose that $\overline\beta(f;\mathcal A)=1$;
then, since $0\not\in\overline\alpha(f;\mathcal A)$, we must have
$|\overline\alpha(f;\mathcal A)|=1$; therefore, upon letting
$f'=f-\lambda$ with $\{\lambda\}=\overline\alpha(f;\mathcal A)$, by (6.4) we
see that $f'$ is a multihyperbolic line which is not a uniline;
clearly $J(f,g)=\nz\Rightarrow J(f',g)=\nz$ which is a contradiction by
(8.2). Thus we must have $\overline\beta(f;\mathcal A)=0$.
Therefore $f$ is a uniline by (5.8.3). Hence $k[f,g]=R_2$ by (8.1), and
$|\alpha(f,g;\mathcal A)|=0$ by (8.3). 

\centerline{}

Here is a two polynomial analogue of (5.9).

\centerline{}

{\bf No Affine Irregular Value Theorem (8.5).} {\it Let $f\in R_2$ be 
$Y$-monic of $Y$-degree $N>0$, and let $g\in R_2$ be such that 
$J(f,g)=\nz$ and 
$|\alpha(f,g;\mathcal A)|=0\not\in\overline\alpha(f;\mathcal A)$.
Then $\beta(f,g;\mathcal A)=0$ and $k[f,g]=R_2$.}

\centerline{}

PROOF. By (5.10) we see that
$|\alpha(f,g;\mathcal A)|=0\not\in\overline\alpha(f;\mathcal A)\Rightarrow
\beta(f,g;\mathcal A)|=0\not\in\overline\alpha(f;\mathcal A)$ and by (8.4)
we see that
$\beta(f,g;\mathcal A)=0\not\in\overline\alpha(f;\mathcal A)\Rightarrow f$
is a uniline.

\centerline{}

Finally here is a two polynomial analogue of (6.4).

\centerline{}

{\bf One Affine Irregular Value Theorem (8.6).} {\it Let $f\in R_2$ be 
$Y$-monic of $Y$-degree $N>0$, and let $g\in R_2$ be such that 
$|\alpha(f,g;\mathcal A)|-1=0\not\in\overline\alpha(f;\mathcal A)$.
Then $J(f,g)\not\in k^{\times}$.}

\centerline{}

PROOF. Let $g'=g-\lambda$ with $\{\lambda\}=\alpha(f,g;\mathcal A)$. Then 
$g'\in R_2$ with 
$\beta(f,g';\mathcal A)=0\not\in\overline\alpha(f;\mathcal A)$
and hence by (8.4) we get $J(f,g')\not\in k^{\times}$. Therefore
$J(f,g)\not\in k^{\times}$.

\centerline{}

{\bf Corollary (8.7).} {\it Let $f\in R_2$ be $Y$-monic of $Y$-degree $N>0$, 
and let $g\in R_2$. Then we have the following.

(8.7.1) If $g\not\in k$ and ${\rm gcd}(f,g-\mu)=1$ for all $\mu\in k$ then
$|\alpha(f,g;\mathcal A)|<\chi_{\infty}(f)$.

(8.7.2) If $0\not\in\overline\alpha(f;\mathcal A)$ and $\chi(F)=2$ 
then $J(f,g)\not\in k^{\times}$.}

\centerline{}

PROOF. In view of (8.5) and (8.6), (8.7.2) follows from (8.7.1) by
noting that  if $g\in k$ then $J(f,g)=0\not\in k^{\times}$, 
and if $f$ and $g-\mu$ for some $\mu\in k$ have a nonconstant 
common factor $\theta$ in $R_2$ 
then $J(f,g)=J(f,g-\mu)$ is divisible by $\theta$ and hence 
$J(f,g)\not\in k^{\times}$. 
To prove (8.7.1), assume that   
$g\not\in k$ and ${\rm gcd}(f,g-\mu)=1$ for all $\mu\in k$.
Then by replacing $g$ by $g-c$ for some $c\in k$ we may suppose
that int$(f,g;\mathcal A)={\rm maxint}(f,g;\mathcal A)$. 
Then int$(f,g;\mathcal A)$ is a
positive integer and, upon letting $(F,G)\sim_m(f,g)$, for all 
$\lambda\in k$ we have 
int$(F,G-\lambda)=-{\rm int}(f,g-\lambda;\mathcal A)$.
Therefore by the inequality for $|\alpha(F,G)|$ given in Section 1 
we conclude that $|\alpha(f,g;\mathcal A)|<\chi_{\infty}(f)$.

\centerline{}

{\bf Remark (8.8).}
Let us note that by a well-known construction, given any finite number
of nonzero elements $f_1,\dots,f_r$ in $R_2$, 
we can find a $k$-automorphism $\sigma$ of $R_2$ such that
$\sigma(f_1),\dots,\sigma(f_r)$ are $Y$-monic.
Clearly it suffices to prove this for their product
$f=f_1\dots f_r$. If $f\in k[Y]$ then we can take $\sigma$ to be identity.
So assume $f\not\in k[Y]$, let $p\ge 0$ be the $(X,Y)$-degree of $f$, 
let $m>0$ be the $X$-degree of $f$, let  $a(Y)$ with $0\ne a(Y)\in k[Y]$ be
the coefficient of $X^m$ in $f$, and let $n\ge 0$ be the $Y$-degree of
$a(Y)$. Take integer $q>p$ and let $\sigma$ be the $k$-automorphism of
$R_2$ given by $\sigma(X)=X+Y^q$ and $\sigma(Y)=Y$. Then
$\sigma(f)$ is $Y$-monic of $Y$-degree $t=n+mq$.
Note that if $r=2$ with $f_i\not\in k[Y]$ for $1\le i\le 2$ and for their
corresponding degrees $p_i,m_i,n_i,t_i$ we have
$m_1/m_2=n_1/n_2$ then we have $t_1/t_2=n_1/n_2$.

\centerline{}

{\bf Remark (8.9).} The AJC = the {\bf Algebraic Jacobian Conjecture}
predicts that for any $f,g$ in $R_2$ with $J(f,g)=\nz$ we have
$k[f,g]=R_2$.
To relate this to the two meromorphic jacobian conjectures, given any
$f,g$ in $R_2$, upon taking $F,G$ in $R$ with $(F,G)\sim_m(f,g)$, and upon 
letting $J=J(F,G)$ and $J\sim_m j'$ with $j'=J(f,g)$, as in the proof 
of (4.6), by the chain rule for jacobians we get $J=-X^{-2}J'$.
To show that Conjecture I implies AJC, in view of the first sentence of the
above Remark (8.8), we may assume that $f\in R_2$ is $Y$-monic of
$Y$-degree $N>0$, and $g\in R_2$ such that $J(f,g)=\nz$, and we want
to show that $k[f,g]=R_2$; now clearly $J=\nz X^{-2}$; by (4.4) and (4.5) we
also have $\beta(f,g;\mathcal A)=\beta(F,G)$ and 
$\overline\alpha(f;\mathcal A)=\overline\alpha(F)$ and hence by Conjecture I
we get $\beta(f,g;\mathcal A)=0\not\in\overline\alpha(f;\mathcal A)$; 
therefore by (8.4) we conclude that $k[f,g]=R_2$.
Before dealing with Conjecture II, let us note that
from what is shown in \cite{Ab3, Ab5} it follows that AJC is equivalent to
a certain variation AJC* of it. To state this, for $1\le i\le 2$ let
$f_i=f_i(X,Y)\in R_2\setminus k[Y]$, let $p_i$ and $m_i$ be the
$(X,Y)$-degree and $X$-degrees of $f_i$, let $0\ne a_i(Y)\in k[Y]$ be the
coefficient of $X^{m_i}$ in $f_i$, let $n_i$ be the $Y$-degree of $a_i(Y)$,
and assume that 
$p_i=m_i+n_i>m_i\text{ and deg}_{(X,Y)}[f_i(X,Y)-a_i(Y)X^{m_i}]< p_i$.
Now AJC* says that, for any such pair $f_1,f_2$, 
if $J(f_1,f_2)=\nz$ and $m_1/m_2=n_1/n_2$ then either
$p_1|p_2$ or $p_2|p_1$. 
[This is a very iffy proposition because, as is
shown in \cite{Ab3, Ab5}, if AJC were true then such a pair $f_1,f_2$
cannot exist]. 
In view of the last sentence of the above Remark (8.8),
to prove AJC* it suffices to prove AJC** which says that if 
$f,g$ in $R_2$ are $Y$-monic of $Y$-degrees $N>0$ and $M>0$ such that
$J(f,g)=\nz$ then either $M|N$ or $N|M$. Clearly Conjecture II implies
AJC**. 

\centerline{}

\centerline{\bf Section 9: Usual Newton Polygon}

\centerline{}

As a tool for dealing with Conjecture II, we shall now revisit the
Usual Newton Polygon as developed in \cite{Ab3}.

Recall that $R=k((X))[Y]$ where $k$ is an algebraically closed field of
characteristic zero, and $R^\natural$ is the set of all irreducible monic 
polynomials of positive degree in $Y$ over $k((X))$.
As in Section 1, for any $F=F(X,Y)\in R$ of $Y$-degree $N$
and branch number $\chi(F)$ we have
$$
F=F_0\prod_{1\le j\le\chi(F)}F_j
$$
where
$$
F_0=F_0(X)\in k((X))
\quad\text{ and }\quad
F_j=F_j(X,Y)\in R^{\natural}\text{ with deg}_YF_j=N_j
$$ 
and for any integer $\nu>0$ which is divisible by 
$N_1,\dots,N_{\chi(F)}$ as Newton factorization of $F$ we have
$$
F(X^{\nu},Y)=F_0(X^\nu)\prod_{1\le i\le N}(Y-z_i(X))
\quad\text{ with }\quad z_i(X)\in k((X)).
$$ 

Note that 
$$
\text{int}(F,Y)=\text{ord}_XF(X,0) 
$$
and let us define the {\bf final root order} of $F$ by putting
$$
\widehat O(F)=\text{max}\{\text{ord}_X z_i(X):1\le i\le N\}
$$
with the understanding that if $N\le 0$ then $O(F)=-\infty$.
For any $c$ which is a rational number or $\infty$ we define the
{\bf vertical label} of $F$ at $c$ and the {\bf starred vertical label} 
of $F$ at $c$ 
to be the nonnegative integers $L(F,c)$ and $L^*(F,c)$ obtained 
by putting
$$
L(F,c)=|\{i\in\{1,\dots,N\}:\text{ord}_X z(i)\ge c\nu\}|
$$
and
$$
L^*(F,c)=\begin{cases}
|\{i\in\{1,\dots,N\}:\text{ord}_X z(i)> c\nu\}|&\text{if }c\ne\infty\\
|\{i\in\{1,\dots,N\}:\text{ord}_X z(i)\ge c\nu\}|&\text{if }c=\infty\\
\end{cases}
$$
with the understanding that if $N\le 0$ then $L(F,c)=0=L^*(F,c)$.
We define the
{\bf final vertical label} and the {\bf postfinal vertical label} of $F$
by putting
$$
\widehat L(F)=L(F,\widehat O(F))
\quad\text{ and }\quad
\widetilde L(F)=L^*(F,\widehat O(F)).
$$
For any integer $a$ we denote the {\bf coefficient} of $X^a$ in $F$ by
coef$_X(F,a)$, i.e., taking summation over all integers $a$ we have
$$
F(X,Y)=\sum \text{coef}_X(F,a)X^a
\quad\text{ with }\quad \text{coef}_X(F,a)\in k[Y]
$$
and we extend this by putting
$$
\text{coef}_X(F,a)=0\text{ if $a$ is either $\infty$ or a rational number
which is not an integer.}
$$
Note that the {\bf $X$-initial coefficient} and the {\bf $X$-initial form}
of $F$ are given by
$$
\text{inco}_XF=\text{coef}_X(F,\text{ord}_X F)
$$
and
$$
\text{info}_X F=\begin{cases}
(\text{inco}_XF)X^a\text{ with }a=\text{ord}_XF&\text{if }F\ne 0\\
0&\text{if }F=0.
\end{cases}
$$
Also note that, given any $y=y(X)\in k((X))$, for all $a$ which is a 
rational number or $\infty$ we have
$$
\text{coef}_X(y,a)\in k.
$$

For any $G=G(X,Y)\in R$ of $Y$-degree $M$
and branch number $\chi(G)$ we have
$$
G=G_0\prod_{1\le l\le\chi(G)}G_l
$$
where
$$
G_0=G_0(X)\in k((X))
\quad\text{ and }\quad
G_l=G_l(X,Y)\in R^{\natural}\text{ with deg}_YG_l=M_l
$$ 
and assuming the integer $\nu>0$ to be divisible by 
$N_1,\dots,N_{\chi(F)},M_1,\dots,M_{\chi(G)}$, in addition to the 
Newton factorization of $F$, as Newton Factorization of $G$ we have
$$
G(X^{\nu},Y)=G_0(X^\nu)\prod_{1\le e\le M}(Y-w_e(X))\quad\text{ with }\quad
w_e(X)\in k((X))
$$ 

We define the {\bf normalized contact} of $F$ and $G$ by putting
$$
\text{noc}(F,G)=(1/\nu)\text{max}\{\text{ord}_X(z_i(X)-w_e(X)):
1\le i\le N\text{ and }1\le e\le M\}
$$
and we note that this is a rational number or $\pm\infty$ and: 
it is $-\infty\Leftrightarrow N\le 0$ or $M\le 0$,
and: it is $\infty\Leftrightarrow F_j=G_l$ for some $j$ and $l$ with
$1\le j\le\chi(F)$ and $1\le l\le\chi(G)$.
We also define the {\bf restricted normalized contact} 
of $F$ and $G$ by putting
\begin{align*}
\text{rnoc}(F,G)=(1/\nu)\text{max}\{\text{ord}_X(z_i(X)-w_e(X))
:1&\le i\le N\text{ and }1\le e\le M\\
&\qquad\text{ with }z_i(X)\ne w_e(X)\}
\end{align*}
and we note that this is a rational number or $-\infty$ and: it is 
$-\infty\Leftrightarrow F_1=\dots=F_{\chi(F)}=G_1=\dots=G_{\chi(G)}$ and
$N_1=\dots=N_{\chi(F)}=1=M_1=\dots=M_{\chi(G)}$.

Assuming $N>0$, to enlarge the pair $\widehat O(F),\widehat L(F)$ 
into the Usual Newton Polygon of $F$, we arrange the set 
$\{(1/\nu)\text{ord}_X z_i(X):1\le i\le N\}$ as an increasing sequence
\begin{align*}
O_1(F)<O_2(F)<\dots&<O_{\iota(F)}(F)\\
&\text{ with preaugumentation } O_0(F)=\text{ord}_X F_0(X)
\end{align*}
and we call $O_i(F)$ the {\bf $i$-th root order} of $F$ 
and $\iota(F)$ the {\bf index} of $F$; note that $\iota(F)$ is the size
of the above set, and $O_0(F),O_1(F),O_2(F),\dots,O_{\iota(F)}(F)$ are 
integers with the exception that $O_{\iota(F)}(F)$ may be $\infty$.
Next we introduce the decreasing sequence of nonnegative integers
\begin{align*}
L_1(F)>L_2(F)>&\dots>L_{\iota(F)}(F)\ge L_{\iota(F)+1}(F)\\
&\quad\text{ with }L_i(F)=L(F,O_i(F))\text{ for }1\le i\le\iota(F)\\
&\quad\text{ and }L_{\iota(F)+1}(F)=\widetilde L(F)
\end{align*}
where we call $L_i(F)$ the {\bf $i$-th vertical label} of $F$. 
The above two sequences together constitute the UNP$(F)=$ the
{\bf Usual Newton Polygon} of $F$. 

To relate the UNP with the customary picture in the real plane,
continuing with the assumption of $N>0$, for any $c$ which ia a
rational number or $\infty$, upon letting
\begin{align*}
(<c)&=\{1\le i\le N:\text{ord}_X z_i(X)<c\nu\},\\
(=c)&=\{1\le i\le N:\text{ord}_X z_i(X)=c\nu\},\\
(>c)&=\{1\le i\le N:\text{ord}_X z_i(X)>c\nu\},
\end{align*}
we define the {\bf horizontal level} $\Lambda(F,c)$ 
of $F$ at $c$ 
and the {\bf starred horizontal level} $\Lambda(F,c)$ 
of $F$ at $c$ by putting
$$
\Lambda(F,c)=O_0(F)
+\left[\sum_{i\in(<c)}(1/\nu)\text{ord}_X z_i(X)\right]
+c|(=c)|
$$
and
$$
\Lambda^*(F,c)=O_0(F)
+\left[\sum_{i\in(<c)}(1/\nu)\text{ord}_X z_i(X)\right]
+c|(=c)|+c|(>c)|
$$
with the understanding that $0$ times $\infty$ is $0$,
and we define the 
{\bf polynomial} $0\ne P^{(F,c)}=P^{(F,c)}(Y)\in k[Y]$ 
of $F$ at $c$ by putting
$$
P^{(F,c)}(Y)=\text{inco}_XF_0(X)
\left[\prod_{i\in(<c)}\text{inco}_X z_i(X)\right]
\left[\prod_{i\in(=c)}(Y-\text{inco}_X z_i(X))\right]
Y^{|(>c)|}.
$$
We define the {\bf final horizontal level} and the {\bf postfinal horizontal
level} of $F$ by putting
$$
\widehat\Lambda(F)=\begin{cases}
\Lambda(F,O_{\iota(F)-1}(F))&\text{if }\iota(F)\ne 1\\
O_0(F)&\text{if }\iota(F)=1\\
\end{cases}
$$
and 
$$
\widetilde\Lambda(F)=\Lambda(F,\widehat O(F))
$$
and we define {\bf final polynomial} of $F$ by putting
$$
\widehat P^{(F)}=\widehat P^{(F)}(Y)=P^{(F,\widehat O(F))}(Y).
$$
We introduce the sequence
$$
\Lambda_1(F)=O_0(F)\quad\text{ and }\quad
\Lambda_i(F)=\Lambda(F,O_{i-1}(F))\text{ for }2\le i\le\iota(F)+1
$$
where we call $\Lambda_i(F)$ the {\bf $i$-th horizontal label} of $F$,
and we introduce the sequence
$$
P_i^{(F)}=P^{(F,O_i(F))}\text{ for }1\le i\le\iota(F)
$$
where we call $0\ne P_i^{(F)}=P_i^{(F)}(Y)\in k[Y]$ the
{\bf $i$-th polynomial} of $F$.

Now the CNP$(F)=$ the {\bf Customary Newton Polygon} of $F$ consists of
the $\iota(F)$ segments in the real plane where, for $1\le i\le\iota(F)$, 
the {\bf $i$-th segment or side} of CNP$(F)$ joins the point 
$(\Lambda_i(F),L_i(F))$ to the point $(\Lambda_{i+1}(F),L_{i+1}(F))$,
with the understanding that if $\widehat O(F)=\infty$
then the last segment or side is the
half-infinite horizontal line emanating from
the point $(\widehat\Lambda(F),\widehat L(F))$ and going to infinity
on the right.  For $1\le i\le\iota(F)$ we embellish the $i$-th side of 
CNP$(F)$ by {\bf placing} the $i$-th polynomial 
$0\ne P_i^{(F)}=P_i^{(F)}(Y)\in k[Y]$ of $F$ on it.
Alternatively, CNP(F) may be constructed thus.
Its first {\bf vertex} is $(\Lambda_1(F),L_1(F))=(O_0(F),N)$. 
The first side is the line of {\bf slope} $O_1(F)$ starting at the
first vertex and ending at height $L_2(F)$ giving us the second
vertex. Inductively, the $i$-the side is defined to be 
the line of slope $O_i(F)$
starting at the $i$-th vertex $(\Lambda_i(F),L_i(F))$ and ending at
height $L_{i+1}(F)$ giving us the $(i+1)$-th vertex whose horizontal
coordinate is defined to be $\Lambda_{i+1}(F)$.
Letting this side continue to height zero, the horizontal
coordinate of the point so obtained is $\Lambda^*(F,O_i(F))$. 
Note that we are interpreting the slope of a side 
to be the tangent of the angle it makes with the $Y$-axis.
{\it This heuristic-geometric paragraph is {\bf not} a logical part
of the paper.} 

Continuing with the assumption of $N>0$, clearly we have
\begin{equation*}
\begin{cases}
\widehat O(F)=O_{\iota(F)}(F)
\text{ and }\widehat P^{(F)}=P^{(F)}_{\iota(F)}P^{(F,\widehat O(F))}\\
\text{ and deg}_Y\widehat P^{(F)}=\widehat L(F)=L_{\iota(F)}(F)
\text{ with }L_1(F)=N
\end{cases}
\tag{9.1}
\end{equation*}
and
\begin{equation*}
\widetilde\Lambda(F)=\Lambda_{\iota(F)+1}(F)
=\text{ord}_X F(X,0)=\text{int}(F,Y)
\text{ with }
\Lambda_1(F)=O_0(F)
\tag{9.2} 
\end{equation*}
and
\begin{equation*}
\widehat\Lambda (F)=\Lambda_{\iota(F)}(F)
\quad\text{ and }\quad
\text{ord}_Y\widehat P^{(F)}=\widetilde L(F)=L_{\iota(F)+1}(F)
\tag{9.3} 
\end{equation*}
and for $1\le i\le\iota(F)$ we have
\begin{equation*}
\aligned
\Lambda_{i+1}(F)
&=\Lambda_i(F)+(L_i(F)-L_{i+1}(F))O_i(F)\\
&=O_0(F)+\sum_{1\le j\le i}(L_j(F)-L_{j+1}(F))O_i(F)
\endaligned
\tag{9.4} 
\end{equation*}
and
\begin{align*}
\text{deg}_Y P_i^{(F)}=L_i(F)
\quad\text{ and }\quad
\text{ord}_Y P_i^{(F)}=L_{i+1}(F)
\tag{9.5} 
\end{align*}
and for any $c$ which is a rational number or $\infty$ we have
\begin{equation*}
\Lambda(F,c)=\begin{cases}
\Lambda_1(F)&\text{if }c<O_1(F)\\
\Lambda_{i+1}(F)&\text{if }O_i(F)\le c<O_{i+1}(F)
\text{ with }1\le i<\iota(F)\\
\Lambda_{\iota(F)+1}(F)&\text{if }O_{\iota(F)}(F)\le c
\end{cases}
\tag{9.6} 
\end{equation*}
and
\begin{equation*}
\Lambda^*(F,c)=\Lambda(F,c)+c|(>c)|
\text{ with $(>c)$ as above}
\tag{9.7} 
\end{equation*}
and
\begin{equation*}
\text{deg}_Y P^{(F,c)}=L(F,c)
\quad\text{ and }\quad
\text{ord}_Y P^{(F,c)}=L^*(F,c).
\tag{9.8} 
\end{equation*}
Moreover, since ord is additive and inco is multiplicative,
for any rational number $c$ for which $c\nu$ is an integer,
we have
\begin{equation*}
\text{ord}_X F(X^\nu,YX^{c\nu})=\Lambda^*(F,c)\nu
\text{ with }
\text{inco}_X F(X^\nu,YX^{c\nu})=P^{(F,c)}(Y).
\tag{9.9} 
\end{equation*}

Assuming $N>0$ and $M>0$,
for $0\le j\le\text{min}(\iota(F),\iota(G))$ we say that 
UNP$(F)$ and UNP$(G)$ are {\bf $j$-step parallel}, 
in symbols we write $\text{UNP}(F)||_j\text{UNP}(G)$, if
$$
\begin{cases}
MO_0(F)=NO_0(G),\text{ and for }1\le i\le j\text{ we have}\\
O_i(F)=O_i(G)\text{ and }ML_i(F)=NL_i(G).
\end{cases}
$$
Moreover, we say that UNP$(F)$ and UNP$(G)$ are {\bf parallel}, 
in symbols we write $\text{UNP}(F)||\text{UNP}(G)$, if
$$
\iota(F)=\iota(G)\text{ and }
\text{UNP}(F)||_{\iota(F)}\text{UNP}(G).
$$
Likewise, we say that UNP$(F)$ is {\bf smaller} than UNP$(G)$, 
in symbols we write $\text{UNP}(F)<\text{UNP}(G)$, if
$$
\begin{cases}
\widehat O(F)<\widehat O(G)\text{ with }\widehat L(G)=1,\text{ and}\\
\text{either }\iota(F)=\iota(G)\text{ with } 
\text{UNP}(F)||_{\iota(F)-1}\text{UNP}(G)
\text{ and }M\widehat L(F)=N\widehat L(G),\\
\,\;\quad\text{or }\iota(F)=\iota(G)-1\text{ with }
\text{UNP}(F)||_{\iota(F)}\text{UNP}(G)
\end{cases}
$$
and we note that
\begin{equation*}
\text{UNP}(F)<\text{UNP}(G)\Rightarrow L^*(G,\widehat O(F))=1.
\tag{*}
\end{equation*}
Finally, we say that UNP$(F)$ and UNP$(G)$ are {\bf pseudoparallel},
in symbols we write $\text{UNP}(F)|.|\text{UNP}(G)$, if
either $\text{UNP}(F)||\text{UNP}(G)$ or $\text{UNP}(F)<\text{UNP}(G)$ or
$\text{UNP}(G)<\text{UNP}(F)$; 
note that these three conditions are mutually exclusive.
In view of (9.1), (9.2) and (9.4), by the definition of parallelness
we see that
\begin{equation*}
\begin{cases}
\text{if $\text{UNP}(F)||\text{UNP}(G)$ then }\\
(M)\text{int}(F,Y)=(N)\text{int}(G,Y)
\text{ and }M\widehat L(F)=N\widehat L(G).
\end{cases}
\tag{9.10}
\end{equation*}

{\bf Calculation.}
Continuing with the assumption that $N>0$ and $M>0$,
{\it let $c$ be a rational number such that $c\nu$ is an integer.}
Let
\begin{equation*}
\widetilde F=\widetilde F(X,Y)=F(X^\nu,YX^{c\nu})
\quad\text{ and }\quad
\widetilde G=\widetilde G(X,Y)=G(X^\nu,YX^{c\nu})
\tag{9.11} 
\end{equation*}
and similarly let
\begin{equation*}
\begin{cases}
\widetilde J=\widetilde J(X,Y)\text{ be obtained by substituting}\\
\qquad\text{$(X^\nu,YX^{c\nu})$ for $(X,Y)$ in $J(F,G)$.}
\end{cases}
\tag{9.12} 
\end{equation*}
Clearly $J(X^\nu,YX^{c\nu})=\nu X^{c\nu+\nu-1}$
and hence by the chain rule for jacobians we get
\begin{equation*}
J(\widetilde F,\widetilde G)=\nu X^{c\nu+\nu-1}\widetilde J(X,Y).
\tag{9.13} 
\end{equation*}
Now
\begin{equation*}
\widetilde F(X,Y)=X^a P(Y)+(\text{terms of $X$-degree}>a)
\tag{9.14}
\end{equation*}
where
\begin{equation*}
a=\text{ord}_X F(X^\nu,YX^{c\nu})
\text{ with }
0\ne P=P(Y)=\text{inco}_X F(X^\nu,YX^{c\nu})\in k[Y]
\tag{9.15}
\end{equation*}
and
\begin{equation*}
\widetilde G(X,Y)=X^b Q(Y)+(\text{terms of $X$-degree}>b)
\tag{9.16}
\end{equation*}
where
\begin{equation*}
b=\text{ord}_X G(X^\nu,YX^{c\nu})
\text{ with }
0\ne Q=Q(Y)=\text{inco}_X G(X^\nu,YX^{c\nu})\in k[Y].
\tag{9.17}
\end{equation*}
Letting $'$ {\bf denote $Y$-derivative} we have
\begin{equation*}
J(\widetilde F,\widetilde G)=X^{a+b-1}(aPQ'-bQP') 
+(\text{terms of $X$-degree}>a+b-1)
\tag{9.18}
\end{equation*}
and hence 
\begin{equation*}
\begin{cases}
\text{ord}_X J(\widetilde F,\widetilde G)\ge a+b-1,\\
\text{with: ord}_X J(\widetilde F,\widetilde G)>a+b-1
\Leftrightarrow aPQ'-bQP'=0
\end{cases}
\tag{9.19}
\end{equation*}
and by (9.13) we see that 
\begin{equation*}
\text{if }J(F,G)\in k((X))
\text{ then }aPQ'-bQP'\in k.
\tag{9.20}
\end{equation*}

\centerline{}

{\bf Related Polynomials.}
To analyze (9.19) and (19.20), let us recall the concept of related 
polynomials developed in \cite{Ab3}. Given any
$$
0\ne P=P(Y)\in k[Y]
\quad\text{ and }\quad
0\ne Q=Q(Y)\in k[Y]
$$
with deg$_YP=n$ and deg$_YQ=m$, we say that $P$ and $Q$ are
{\bf related} to mean that $P^m=\nz Q^n$.
Recall that upon letting
$$
P=P(Y)=P_0\prod_{1\le i\le p}(Y-u_i)^{r_i}
\quad\text{ and }\quad
Q=Q(Y)=Q_0\prod_{1\le j\le q}(Y-v_j)^{s_j}
$$
with $P_0\ne 0\ne Q_0$ in $k$, 
pairwise distinct elements  $u_1,\dots,u_p$ in $k$, 
pairwise distinct elements  $v_1,\dots,v_q$ in $k$, 
nonnegative integers $p,q$, and positive 
integers $r_1,\dots,r_p,s_1,\dots,s_q$, 
we have
$$
\text{rad}(P)=P_0\prod_{1\le i\le p}(Y-u_i)
\quad\text{ and }\quad
\text{rad}(Q)=Q_0\prod_{1\le j\le q}(Y-v_j)
$$
and note that
\begin{align*}
&PQ\text{ has a multiple root}\\
&\Leftrightarrow 
\text{ either $r_i\ge 2$ for some $i$ or $s_j\ge 2$
for some $j$ or $u_i=v_j$ for some $i,j$.}
\end{align*}
Clearly
\begin{equation*}
\text{if $m+n\ne 0$ $P$ and $Q$ are related then $n\ne 0\ne m$}
\tag{9.21}
\end{equation*}
and by a standard argument we see that
\begin{equation*}
\begin{cases}
\text{if $m+n\ne 0$ then: $P$ and $Q$ are related }
\Leftrightarrow\\
\text{rad}(P)=\nz\text{rad}(Q)
\text{ and upon relabelling $v_1,\dots,v_q$}\\
\text{so that $u_i=v_i$ for $1\le i\le p=q$}\\
\text{we have }mu_i=nv_i\text{ for }1\le i\le p
\end{cases}
\tag{9.22}
\end{equation*}
and
\begin{equation*}
\begin{cases}
\text{if $m+n\ne 0$ and $a,b$ are integers}\\
\text{such that $aPQ'-bQP'=0$ with either $a<0$ or $b<0$}\\
\text{then $P^{-b}=\nz Q^{-a}$ with $a<0$ and $b<0$}\\
\text{and $P$ and $Q$ are related with $ma=nb$.}
\end{cases}
\tag{9.23}
\end{equation*}
Clearly
\begin{equation*}
\begin{cases}
\text{if $a,b$ are integers}\\ 
\text{such that $aPQ'-bQP'\in k^{\times}$ 
with either $a<0$ or $b<0$}\\
\text{then $PQ$ has no multiple root}
\end{cases}
\tag{9.24}
\end{equation*}
and hence by (9.23) we see that
\begin{equation*}
\begin{cases}
\text{if $PQ$ has a multiple root and $a,b$ are integers}\\
\text{such that $aPQ'-bQP'\in k$ with either $a<0$ or $b<0$}\\
\text{then $aPQ'-bQP'=0$}\\
\text{and $P^{-b}=\nz Q^{-a}$ with $a<0$ and $b<0$}\\
\text{and $P$ and $Q$ are related with $ma=nb$.}
\end{cases}
\tag{9.25}
\end{equation*}

{\bf Calculation Continued.} 
Reverting to the definition of $a,b,P,Q$ given in (9.15) 
and (9.17), we continue to let deg$_YP=n$ and deg$_YG=m$.
By (9.2), (9.4), (9.6), (9.7) and  (9.9) we see that
\begin{equation*}
a\le\text{int}(F,Y)\nu
\quad\text{ and }\quad
b\le\text{int}(G,Y)\nu
\tag{9.26}
\end{equation*}
and hence
\begin{equation*}
\begin{cases}
\text{if either int$(F,Y)<0$ or int$(G,Y)<0$}\\ 
\text{then either $a<0$ or $b<0$.}\\ 
\end{cases}
\tag{9.27}
\end{equation*}
By (9.8) we see that
\begin{equation*}
\begin{cases}
\text{(i) deg}_YP>0
\Leftrightarrow\text{deg}_YP\ge\widehat L(F)
\Leftrightarrow c\le\widehat O(F),\\
\text{(ii) ord}_YP>0
\Leftrightarrow\text{ord}_YP\ge\widehat L(F)
\Leftrightarrow c<\widehat O(F),\\
\text{(ii) deg}_YP>\text{ord}_Y P\Leftrightarrow 
c\in\{O_1(F),\dots,O_{\iota(F)}(F)\},\\
\text{(iv) deg}_YQ>0
\Leftrightarrow\text{deg}_YQ\ge\widehat L(G)
\Leftrightarrow c\le\widehat O(G),\\
\text{(v) ord}_YQ>0
\Leftrightarrow\text{ord}_YQ\ge\widehat L(G)
\Leftrightarrow c<\widehat O(G),\\
\text{(vi) deg}_YQ>\text{ord}_Y Q\Leftrightarrow 
c\in\{O_1(G),\dots,O_{\iota(G)}(G)\},\\
\text{(vii) }c< O_{\iota(F)-1}\text{ with }\iota(F)\ge 2
\Rightarrow\text{ord}_YP\ge 2,\\
\text{(viii) }c< O_{\iota(G)-1}\text{ with }\iota(G)\ge 2
\Rightarrow\text{ord}_YQ\ge 2,\\
\end{cases}
\tag{9.28}
\end{equation*}
and hence
\begin{equation*}
\begin{cases}
\text{if either (i) }c<\text{min}(\widehat O(F),\widehat O(G)),\\
\text{or (ii) }c<\widehat O(G))\text{ with }\widehat L(G)\ge 2,\\
\text{or (iii) }c<O_{\iota(G)-1}(G))\text{ with }\iota(G)\ge 2,\\
\text{or (iv) }c<\widehat O(F))\text{ with }\widehat L(F)\ge 2,\\
\text{or (v) }c<O_{\iota(F)-1}(F))\text{ with }\iota(F)\ge 2,\\
\text{then ord}_Y PQ\ge 2,
\end{cases}
\tag{9.29}
\end{equation*}
and
\begin{equation*}
\begin{cases}
\text{if either }
c\in\{O_1(F),\dots,O_{\iota(F)}(F)\}
\setminus \{O_1(G),\dots,O_{\iota(G)}(G)\}\\
\text{or }c\in\{O_1(G),\dots,O_{\iota(G)}(G)\}
\setminus \{O_1(F),\dots,O_{\iota(F)}(F)\}\\
\text{then $P$ and $Q$ are not related.}
\end{cases}
\tag{9.30}
\end{equation*}

\centerline{}

{\bf Main Lemma (9.31).} {\it Let $F$ and $G$ in $R$ be 
of $Y$-degrees $N>0$ and $M>0$ respectively,
and assume that $J(F,G)\in k((X))$ and 
${\rm UNP}(F)||_0{\rm UNP}(G)$. Also assume that
either {\rm int}$(F,Y)<0$ or {\rm int}$(G,Y)<0$. Then
${\rm UNP}(F)|.|{\rm UNP}(G)$, and for
$1\le i<{\rm min}(\iota(F),\iota(G))$ the polynomials
$P_i^{(F)}$ and $P_i^{(G)}$ are related. Moreover, if
$\widehat O(F)=\widehat O(G)$ then ${\rm UNP}(F)||{\rm UNP}(G)$
and hence in particular $(M){\rm int}(F,Y)=(N){\rm int}(G,Y)$
and $M\widehat L(F)=N\widehat L(G)$.}

\centerline{}

PROOF. By induction on $j$ we shall show that, given any integer
$j$ with $0\le j\le\text{min}(\iota(F),\iota(G))$, we have
$(1_j)$ to $(6_j)$ stated below. By taking
$j=\text{min}(\iota(F),\iota(G))$ this will establish the lemma.

$(1_j)$ If $j<\text{min}(\iota(F),\iota(G))$ then we have:
$\text{UNP}(F)||_j\text{UNP}(G)$, and the polynomials
$P_i^{(F)}$ and $P_i^{(G)}$ are related for $1\le i\le j$. 

$(2_j)$ If $\widehat O(F)=\widehat O(G)$ and
$j=\text{min}(\iota(F),\iota(G))$ then we have:
$j=\iota(F)=\iota(G)$,
$\text{UNP}(F)||_j\text{UNP}(G)$, the polynomials
$P_i^{(F)}$ and $P_i^{(G)}$ are related for $1\le i<j$, 
and $(M)\text{int}(F,Y)=(N)\text{int}(G,Y)$.

$(3_j)$ If $\widehat O(F)<\widehat O(G)$ and
$j=\text{min}(\iota(F),\iota(G))=\iota(G)$ then we have:
$j=\iota(F)$, $\text{UNP}(F)||_{j-1}\text{UNP}(G)$, 
the polynomials
$P_i^{(F)}$ and $P_i^{(G)}$ are related for $1\le i<j$, 
$M\widehat L(F)=N\widehat L(G)$,
and $\widehat L(G)=1$.

$(4_j)$ If $\widehat O(F)<\widehat O(G)$ and
$j=\text{min}(\iota(F),\iota(G))\ne\iota(G)$ then we have:
$j=\iota(F)=\iota(G)-1$, $\text{UNP}(F)||_j\text{UNP}(G)$, 
the polynomials
$P_i^{(F)}$ and $P_i^{(G)}$ are related for $1\le i<j$, 
and $\widehat L(G)=1$.

$(5_j)$ If $\widehat O(G)<\widehat O(F)$ and
$j=\text{min}(\iota(F),\iota(G))=\iota(F)$ then we have:
$j=\iota(G)$, $\text{UNP}(F)||_{j-1}\text{UNP}(G)$, 
the polynomials
$P_i^{(F)}$ and $P_i^{(G)}$ are related for $1\le i<j$, 
$M\widehat L(F)=N\widehat L(G)$,
and $\widehat L(F)=1$.

$(6_j)$ If $\widehat O(G)<\widehat O(F)$ and
$j=\text{min}(\iota(F),\iota(G))\ne\iota(F)$ then we have:
$j=\iota(G)=\iota(F)-1$, $\text{UNP}(F)||_j\text{UNP}(G)$, 
the polynomials
$P_i^{(F)}$ and $P_i^{(G)}$ are related for $1\le i<j$, 
and $\widehat L(F)=1$.

By hypothesis 
this holds for $j=0$. 
So let $j>0$ and assume for $j-1$. 
Note that now  $(2_{j-1})$ to $(6_{j-1})$ are vacuous, and so
in proving $(1_j)$ to $(6_j)$ we shall only use $(1_{j-1})$ 
and that we shall do without mentioning it explicitly;
we shall also tacitly use the fact that $ML_j(F)=NL_j(G)$ 
which in case of $j>1$ follows from (9.5) and the relatedness of 
$P^{(F)}_{j-1}$ and $P^{(G)}_{j-1}$, and is obvious in case of $j=1$
because $L_1(F)=N$ and $L_1(G)=M$.
Since either int$(F,Y)<0$ or int$(G,Y)<0$, upon letting
$c=\text{min}(O_j(F),O_j(G))$ we see that $c$ is a rational number
such that $c\nu$ is an integer. So we may use the above
Calculation, and then by (9.20) we see that
$aPQ'-bQP'\in k$ and by (9.27) we see that 
$a<0$ and $b<0$.
If $j<\text{min}(\iota(F),\iota(G))$ then by (9.29)(i) we see that
$PQ$ has a multiple root and therefore by
(9.25) and (9.30) we see that $P=P_j^{(F)}$ and $Q=P_j^{(G)}$ are 
related with $O_j(F)=O_j(G)$; this proves $(1_j)$.
If $\widehat O(F)=\widehat O(G)$ and
$j=\text{min}(\iota(F),\iota(G))$ then obviously
$j=\iota(F)=\iota(G)$ with 
$O_j(F)=\widehat O(F)=\widehat O(G)=O_j(G)=\widehat O(G)$
and hence by (9.10) we see that
$(M)\text{int}(F,Y)=(N)\text{int}(G,Y)$; this proves $(2_j)$.
If $\widehat O(F)<\widehat O(G)$ and
$j=\text{min}(\iota(F),\iota(G))=\iota(G)$ then
by (9.30) we see that $P$ and $Q$ are not related and hence by
(9.25) and (9.29)(i,ii) we get $j=\iota(F)$ and
$\widehat L(G)=1$; this proves $(3_j)$.
If $\widehat O(F)<\widehat O(G)$ and
$j=\text{min}(\iota(F),\iota(G))\ne\iota(G)$ then
$j=\iota(F)<\iota(G)$ and by (9.28)(ii,v) we see that
$P$ and $Q$ are not related and hence by (9.25) and
(9.29)(ii,iii) we get $j=\iota(F)=\iota(G)-1$
with $O_j(F)=O_j(G)$ and $\widehat L(G)=1$; this proves $(4_j)$.
Interchanging $F$ and $G$ in the proof of $(3_j)$ and $(4_j)$
we get $(5_j)$ and $(6_j)$.

\centerline{}

{\bf Definition (9.32).}
Using (9.31), in (9.33) to (9.38) we shall show that, under certain 
condition, $J(F,G)\in k((X))$ implies that most branches of 
$F$ and $G$ can be partitioned into {\bf packets} 
$(F_1,\dots,F_r,G_1,\dots,G_s)$ whose members are {\bf pseudocognates}
of each other, i.e., their roots coincide 
upto the last characteristic term. To define 
these concepts more precisely, let us review some relevant terms.

For any $f\in R^{\natural}$ of $Y$-degree $n$, the 
{\bf newtonian sequence of characteristic exponents} 
of $f$ relative to $n$, denoted by 
$$
m(f)=m_i(f)_{0\le i\le h(m(f))+1}
$$
is defined on pages 
page 3-4 of \cite{AA2},
where the {\bf GCD-sequence} 
$$
d(m(f))=d_i(m(f))_{0\le i\le h(d(m(f)))+2}
$$ 
of $m(f)$ is also defined.
For simplicity of notation we put 
$$
h(f)=h(m(f))=h(m(d(f)))
\quad\text{ and }\quad
d(f)=d(m(f))
$$
and
$$
\widehat d(f)=d_{h(f)}(f).
\quad\text{ and }\quad
\widehat m(f)=\begin{cases}
m_{h(f)}(f)&\text{if }h(f)\ne 0\\
-\infty&\text{if }h(f)=0\\
\end{cases}
$$
Note that then 
$$
\quad d_{h(f)+1}(f)=1
\quad\text{ and }
m_{h(f)+1}(f)=\infty
$$ 
and
$$
d_0(f)=0
\quad\text{ and }\quad
m_0(f)=n=d_1(f)
$$ 
and
$$
h(f)=0\Leftrightarrow f(X,Y)=Y\Leftrightarrow\widehat d(f)=0
\Leftrightarrow\widehat m(f)=-\infty.
$$
Also note that:

(9.32.0) rnoc$(f,f)=\widehat m(f)/n$.

For any $c$ which is a rational number or $\infty$ we define the 
{\bf $c$-position} of $f$ to be the unique nonnegative integer $p(f,c)\le h(f)$
such that $m_i(f)/n<c$ for $1\le i\le p(f,c)$, and $c\le m_j(f)/n$ for
$p(f,c)<j\le h(f)$. We also we define the quantities
$$
\widehat d(f,c)=d_{p(f,c)+1}(f)
\quad\text{ and }\quad
\widehat m(f,c)=m_{p(f,c)+1}(f)
$$
and 
$$
t(f,c)=\begin{cases}
\text{the minimal monic polynomial of}\\
\sum_{i<cn}\text{coef}_X(\eta(X),i)X^i\text{ over }k((X^n))\\
\text{were $\eta(X)$ is a root of $f(X^n,Y)$ in $k((X))$}
\end{cases} 
$$
where we call $t(f,c)=t(f,c)(X,Y)\in R^{\natural}$ the 
{\bf $c$-normalized truncation} of $f$. 
Note that then 
$$
\widehat d(f,c)=n/\text{deg}_Y t(f,c)
$$
and
$$
\begin{cases}
\text{if either $h(f)\ne 0$ with $m_{h(f)}(f)/n<c$ or $h(f)=0$,}\\
\text{then $p(f,c)=h(f)$ with $\widehat m(f,c)=\infty$ and $\widehat d(f,c)=1$}
\end{cases}
$$
and
$$
\begin{cases}
\text{if $h(f)\ne 0$ with $m_{h(f)}(f)/n=c$,}\\
\text{then $p(f,c)=h(f)-1$ with $\widehat m(f,c)=\widehat m(f)$ 
and $\widehat d(f,c)=\widehat d(f)$.}
\end{cases}
$$

Let $f'\in R^{\natural}$ be of $Y$-degree $n'$. 
We say that $f$ is a {\bf cognate} of $f'$ if 
noc$(f,f')=\widehat m(f)/n=\widehat m(f')/n'$.
We say that $f$ is an {\bf overcognate} of $f'$ if $h(f)\ne 0\ne h(f')$ with 
noc$(f,f')>\text{max}(\widehat m(f)/n,\widehat m(f')/n')$.
We say that $f$ is a {\bf subcognate} of $f'$ if
$\widehat m(f)/n<\text{noc}(f,f')=\widehat m(f')/n'$.
Note that:

(9.32.1) if $f$ is a cognate (resp: overcognate) of $f'$ then $f'$ is a 
cognate (resp: overcognate) of $f$; 

(9.32.2) if $f$ is a cognate of $f'$ then $h(f)\ne 0\ne h(f')$;

(9.32.3) if $h(f)=0=h(f')$ then $f$ is an overcognate of $f'$;

(9.32.4) if $f$ is an overcognate of $f'$ then 
$\widehat m(f)/n=\widehat m(f')/n'$;

(9.32.5) if $f$ is a cognate or overcognate of $f'$ then $n=n'$ and 
$h(f)=h(f')$ with $m(f)=m(f')$ and $d(f)=d(f')$ and 
$\widehat d(f)=\widehat d(f')$, and for any rational number 
$c$ we have $p(f,c)=p(f',c)$ with $\widehat m(f,c)=\widehat m(f',c)$
and $\widehat d(f,c)=\widehat d(f',c)$;

(9.32.6) if $h(f)+1=h(f')$ with noc$(f,f')=\widehat m(f')/n'$ then
$f$ is a subcognate of $f'$; and

(9.32.7) if $f$ is a subcognate of $f'$ then $h(f)+1=h(f')$ and
$m_i(f)/n=m_i(f')/n'$ for $1\le i\le h(f)$ with
coef$_X(\eta(X),\widehat m(f')n/n')=0$ 
where $\eta(X)$ is a root of $f(X^n,Y)$ in $k((X))$.

Also note that if $h(f')\ne 0$ and $f$ is the $\widehat d(f')$-th approximate 
root of $f'$ in the sense of \cite{Ab2} then $f$ is a subcognate of $f'$;
consequently we may think of a subcognate of $f'$ as a 
{\bf last pseudoapproximate root} of $f'$.
Moreover, if either $f$ is a cognate of $f'$, or $f$ is an overcognate of $f'$,
or $f$ is a subcognate of $f'$, or $f'$ is a subcognate of $f$, then we may
think of $f$ and $f'$ as being {\bf pseudocognates} of each other.

By an {\bf equilateral sequence} in $R^{\natural}$
we mean a sequence $(f_i)_{1\le i\le r}$ of
members of $R^{\natural}$, with integer $r>1$,
for which there exists a (necessarily unique) rational number $c$ such that
for all $i\ne j$ in $\{1,\dots,r\}$ we have noc$(f_i,f_j)=c$ and for all 
$i\in $ $\{1,\dots,r\}$ we have rnoc$(f_i,f_i)\le c$; we call $c$ the
{\bf diameter} of the sequence.
By a {\bf cognate sequence} in $R^{\natural}$ we mean an equilateral sequence 
$(f_i)_{1\le i\le r}$ in $R^{\natural}$
such that for all $i\ne j$ in $\{1,\dots,r\}$ we have that 
$f_i$ is a cognate of $f_j$.
By an {\bf overcognate sequence} in $R^{\natural}$ 
we mean an equilateral sequence 
$(f_i)_{1\le i\le r}$ in $R^{\natural}$
such that for all $i\ne j$ in $\{1,\dots,r\}$ we have that 
$f_i$ is an overcognate of $f_j$.
By a {\bf subcognate sequence} in $R^{\natural}$ 
we mean an equilateral sequence 
$(f_i)_{1\le i\le r}$ in $R^{\natural}$
for which there exists a unique $i'$ in $\{1,\dots,r\}$ such that for all 
$j\ne i'$ in $\{1,\dots,r\}$ we have that $f_{i'}$ is a subcognate of $f_j$
and for all $i\ne j$ in $\{1,\dots,r\}\setminus\{i'\}$ we have that 
$f_i$ is a cognate of $f_j$; we call $f_{i'}$ the {\bf special} branch of the
sequence.
By an {\bf equicognate sequence} in $R^{\natural}$ we mean 
an {\bf equilateral sequence} in $R^{\natural}$ which is either 
a cognate sequence in $R^{\natural}$ or
an overcognate sequence in $R^{\natural}$ or
a subcognate sequence in $R^{\natural}$. 

Note that for an equicognate sequence 
$(f_i)_{1\le i\le r}$ in $R^{\natural}$ with diameter $c$ we have that:

(9.32.8) if the sequence is cognate then for $1\le i\le r$ we have
$$
\widehat d(f_i,c)=\widehat d(f_i)=\widehat d(f_1)
\quad\text{ and }\quad
d_1(f_i)=d_1(f_1);
$$

(9.32.9) if the sequence is overcognate then for $1\le i\le r$ we have
$$
\widehat d(f_i,c)=1
\quad\text{ and }\quad
d_1(f_i)=d_1(f_1);
$$

(9.32.10) and if the sequence is subcognate and we have labelled the branches 
so that the special branch is $f_r$ then for $1\le i< r$ we have
$$
\begin{cases}
1=\widehat d(f_r,c)\le\widehat d(f_i,c)=\widehat d(f_i)=\widehat d(f_1)\\
\text{and }d_1(f_r)=d_1(f_i)/\widehat d(f_i)=d_1(f_1)/\widehat d(f_1).
\end{cases}
$$

Let $F$ and $G$ in $R$ be $Y$-monic of $Y$-degrees $N>0$ and $M>0$ 
respectively. By a {\bf bisequence} of $(F,G)$ we mean a pair of
families $(F_j)_{j\in J},(G_l)_{l\in J'}$ where $J$ and $J'$ are nonempty
subsets of $\{1,\dots,\chi(F)\}$ and $\{1,\dots,\chi(G)\}$ respectively.
This induces the sequence $(f_i)_{1\le i\le r}$ in $R^{\natural}$ where
$r=|J|+|J'|$ and, upon letting $j_1<\dots<j_{|J|}$ and $l_1<\dots<l_{|J'|}$
be the increasing labellings of $J$ and $J'$ respectively, we have
$f_i=F_{j_i}\text{ for }1\le i\le |J|\text{ and }
f_{|J|+i}=G_{l_i}\text{ for }1\le i\le |J'|$.
The bisequence is said to be {\bf equilateral}, ..., {\bf equicognate} if the 
induced sequence is respectively {\bf equilateral}, ..., {\bf equicognate}.
By the {\bf diameter} of an equilateral bisequence we mean the diameter of
the induced sequence. By the {\bf special} branch of a subcognate bisequence
we mean the special branch of the induced sequence.

An equilateral bisequence $(F_j)_{j\in J},(G_l)_{l\in J'}$ of $(F,G)$,
whose diameter is $c$, is said to be {\bf saturated} if: 

(9.32.11) there is no $i\in\{1,\dots,\chi(F)\}\setminus J$
such that either for some $j\in J$ we have noc$(F_i,F_j)\ge c$ or
for some $l\in J'$ we have noc$(F_i,G_l)\ge c$; 

(9.32.12) there is no $i'\in\{1,\dots,\chi(G)\}\setminus J'$
such that either for some $j\in J$ we have noc$(G_{i'},F_j)\ge c$ or
for some $l\in J'$ we have noc$(G_{i'},G_l)\ge c$;

(9.32.13) and for all $j\in J$ and $l\in J'$ we have
noc$(F_j,G)=c$ and noc$(F,G_l)=c$. 
\noindent

An equilateral bisequence $(F_j)_{j\in J},(G_l)_{l\in J'}$ of $(F,G)$,
whose diameter is $c$, is said to be {\bf balanced} if it is saturated and: 

(9.32.14) for all $j\in J$ and $l\in J'$ we have
int$(F_j,G)<0$ and int$(F,G_l)<0$ with
$$
\frac{{\rm int}(F,G_l)}{{\rm int}(F_j,G)}
=\frac{NM_l}{MN_j};
$$

(9.32.15) and we have
$$
\frac{\sum_{j\in J}\widehat d(F_j,c)}{\sum_{l\in J'}\widehat d(G_l,c)}
=\frac{N}{M}.
$$

An equilateral bisequence $(F_j)_{j\in J},(G_l)_{l\in J'}$ of $(F,G)$
is said to be {\bf well-balanced} if it is balanced
and for it:

(9.32.16) the degrees satisfy the equation
$$
\frac{\sum_{j\in J}N_j}{\sum_{l\in J'}M_l}
=\frac{N}{M};
$$

(9.32.17) the intersection multiplicities satisfy the equation
$$
\sum_{j\in J}\text{int}(F_j,G)=\sum_{l\in J'}\text{int}(F,G_l);
$$

(9.32.18) there exist unique negative rational numbers $N'$ and $M'$
with $MN'=NM'$ such that for all $j\in J$ and $l\in J'$ we have
$N_j=N'\text{int}(F_j,G)$ and $M_l=M'\text{int}(F,G_l)$; and 

(9.32.19) upon letting $E=\text{max}((N_j)_{j\in J},(M_l)_{l\in J'})$
and $D=\text{min}(N_j)_{j\in J}$ and $D'=\text{min}(M_l)_{l\in J'}$
we have that: 

(i) if $D\ne E$ then there is a unique $s\in J$ such that
$N_s|E$ with $N_s<E=N_j=M_l$ for all 
$j\in J\setminus\{s\}$ and $l\in J'$, and

(ii) if $D'\ne E$ then there is a unique $s'\in J'$ such that
$M_{s'}|E$ with $M_{s'}<E=N_j=M_l$ for all 
$j\in J$ and $l\in J'\setminus\{s'\}$.

Finally, by a {\bf packet} of $(F,G)$ we mean a 
balanced equicognate bisequence $(F_j)_{j\in J},(G_l)_{l\in J'}$ of $(F,G)$.
Note that then a packet has properties (9.32.11) to (9.32.19), and its
induced sequence has properties (9.32.8) to (9.32.10).
In view of the last 3 references, (i) or (ii) occur exactly when the packet 
is subcognate with $E\ne 1$, and then $F_s$ or $G_{s'}$ is the special
branch. In all other cases, i.e., if the packet is subcognate with $E=1$, 
or the packet is cognate, or the packet is overcognate, then for all
$j\in J$ and $l\in J'$ we have $N_j=E=G_l$.  

\centerline{}

{\bf First Corollary (9.33).}
{\it Let $F$ and $G$ in $R$ be $Y$-monic of $Y$-degrees $N>0$ and $M>0$ 
respectively, and assume that $J(F,G)\in k((X))$. Let $y(X)\in k((X))$ 
be such that ${\rm ord}_XF(X^\nu,y(X))<0$. Let
$$
\begin{cases}
c=(1/\nu){\rm max}_{1\le i\le N}{\rm ord}_X(y(X)-z_i(X)),\\
I=\{i:1\le i\le N\text{ with }{\rm ord}_X(y(X)-z_i(X))=c\nu\},\\
J=\{j:1\le j\le\chi(F)\text{ with }F_j(X^\nu,z_i(X))=0
\text{ for some }i\in I\},\\
\end{cases}
$$
and
$$
\begin{cases}
c'=(1/\nu){\rm max}_{1\le e\le M}{\rm ord}_X(y(X)-w_e(X)),\\
I'=\{e:1\le e\le M\text{ with }{\rm ord}_X(y(X)-w_e(X))=c'\nu\}\\
J'=\{l:1\le l\le\chi(G)\text{ with }G_l(X^\nu,w_e(X))=0
\text{ for some }e\in I'\}.\\
\end{cases}
$$
Then $c$ is a rational number such that $c\nu$ is an integer,
and we have the following.

(9.33.1) If $c<c'$ then $|I(G)|=1$.

(9.33.2) If $c=c'$ then ${\rm ord}_XG(X^\nu,y(X))<0$ and
$$
\frac{{\rm ord}_X F(X^\nu,y(X))}{{\rm ord}_X G(X^\nu,y(X))}
=\frac{N}{M}
=\frac{\sum_{j\in J}\widehat d(F_j,c)}
{\sum_{l\in J'}\widehat d(G_l,c)}.
$$

(9.33.3) If $c=c'$, and there exists $\epsilon\in I'$ 
such that for all $i\in I$ we have
${\rm coef}_X(z_i(X),c\nu)\ne{\rm coef}_X(w_\epsilon(X),c\nu)$,
then
$$
\begin{cases}
\text{for all $i\in I$ and $e\in I'$ we have
${\rm coef}_X(z_i(X),c\nu)\ne{\rm coef}_X(w_e(X),c\nu)$,}\\
\text{for all $i\ne i'$ in $I$ we have
${\rm coef}_X(z_i(X),c\nu)\ne{\rm coef}_X(z_{i'}(X),c\nu)$, and}\\
\text{for all $e\ne e'$ in $I'$ we have
${\rm coef}_X(w_e(X),c\nu)\ne{\rm coef}_X(w_{e'}(X),c\nu)$.}\\
\end{cases}
$$
}

\centerline{}

PROOF.
Let $\overline c=c\nu$ and $\overline c'=c'\nu$. 
Since ${\rm ord}_XF(X^\nu,y(X))<0$, it follows that $c$ is a rational 
number such that $\overline c$ is an integer. Let
$$
\overline X=X^{\nu}
\quad\text{ and }\quad
\overline Y=Y+y(X)
$$
and
$$
\overline F=\overline F(X,Y)=F(\overline X,\overline Y)
\quad\text{ and }\quad
\overline G=\overline G(X,Y)=G(\overline X,\overline Y).
$$
By the chain rule for jacobians we get 
$J(\overline F,\overline G)\in k((X))$. 
Clearly int$(\overline F,Y)={\rm ord}_XF(X^\nu,y(X))$ and hence
int$(\overline F,Y)<0$. 
Moreover $\widehat O(\overline F)=\overline c$ 
and $\widehat O(\overline G)=\overline c'$.
Now by taking $\overline F,\overline G$ for $F,G$ in (9.31) we see that 
if $\widehat O(\overline F)<\widehat O(\overline G)$ 
then $L^*(\overline G,\widehat O(F))=1$. From this it follows that 
ord$_X(y(X)-w_e(X))>\overline c$ for at most one $e$ in $\{1,\dots,M\}$, 
which proves (9.33.1). 

Henceforth assume that $c=c'$.
Again by taking $\overline F,\overline G$ for $F,G$ in (9.31) we see that 
$\text{UNP}(\overline F)||\text{UNP}(\overline G)$ and 
$\text{int}(\overline G,Y)<0$ with
$$
\frac{\text{int}(\overline F,Y)}{\text{int}(\overline G,Y)}
=\frac{N}{M}=\frac{\widehat L(\overline F)}{\widehat L(\overline G)}.
$$
Clearly int$(\overline G,Y)={\rm ord}_XG(X^\nu,y(X))$ and hence we get
the first equality of (9.33.2). Also clearly
$$
\widehat L(\overline F)=|I|=\sum_{j\in J}\widehat d(F_j,c)
\quad\text{ and }\quad 
\widehat L(\overline G)=|I'|=\sum_{l\in J'}\widehat d(G_l,c)
$$
and hence we get the second equality of (9.33.2).
Taking $\overline F,\overline G,\overline c,1$ for $F,G,c,\nu$ 
in the above Calculation we have
$\widetilde F=\widetilde F(X,Y)=\overline F(X,YX^{\overline c})$ and 
$\widetilde G=\widetilde G(X,Y)=\overline G(X,YX^{\overline c})$. 
Clearly the members of the sets
$$
\{\text{coef}_X(z_i(X),\overline c)-\text{coef}_X(y(X),\overline c):i\in I\}
$$
and
$$
\{\text{coef}_X(w_e(X),\overline c)-\text{coef}_X(y(X),\overline c):e\in I'\}
$$
are precisely the roots of $P(Y)$ and $Q(Y)$.
Now assume that there exists $\epsilon\in I'$ 
such that for all $i\in I$ we have
${\rm coef}_X(z_i(X),c\nu)\ne{\rm coef}_X(w_\epsilon(X),c\nu)$,
Then
$\text{coef}_X(w_{\epsilon},\overline c)-\text{coef}_X(y(X),\overline c)$
belongs to the second set
but not to the first set, and hence $P$ and $Q$ are not related.
Therefore (9.33.3) follows from (9.25) and (9.27).

\centerline{}

{\bf Second Corollary (9.34).}
{\it Let $F$ and $G$ in $R$ be $Y$-monic of $Y$-degrees $N>0$ and $M>0$ 
respectively, and assume that $J(F,G)\in k((X))$. For an integer $v$ 
with $1\le v\le\chi(G)$ let $c={\rm noc}(F,G_v)$
and assume that {\rm int}$(F,G_v)<0$. Then $c$ is a rational number such
that $c\nu$ is an integer, and upon letting
$$
J=\{j:1\le j\le\chi(F)\text{ with }{\rm noc}(F_j,G_v)\ge c\}
$$
and
$$
J'=\{j:1\le l\le\chi(G)\text{ with }{\rm noc}(G_l,G_v)\ge c\}
$$
we have that $(F_j)_{j\in J},(G_l)_{l\in J'}$ is a 
balanced equilateral bisequence of $(F,G)$ with diameter $c$. In particular
$J\ne\emptyset$ and $v\in J'$.}

\centerline{}

PROOF. 
Since int$(F,G_v)<0$, we must have $c\ne\infty$. Therefore $c$
is a rational number and upon
letting $\overline c=c\nu$ we see that $\overline c$ is an integer.
Clearly $J\ne\emptyset$ and $v\in J'$. Also clearly, for every $j\in J$
we have noc$(F_j,G_v)=c$.
We can take $\epsilon$ with $1\le \epsilon\le M$ such that
$G_v(X^\nu,w_\epsilon(X))=0$. Now by taking $w_\epsilon(X)$ for $y(X)$ in 
(9.33.1) we see that rnoc$(G_v,G_v)\le c$ and for all $l\ne v$ in $J'$ we
have noc$(G_l,G_v)=c$. 
We can take $\lambda\in k$ such that
for $1\le i\le N$ and $1\le e\le M$ we have 
$$
\text{coef}_X(z_i(X),\overline c)\ne\lambda
\ne\text{coef}_X(w_e(X),\overline c).
$$ 
Henceforth let $y(X)\in k((X))$ be defined by putting
$$
y(X)=\lambda X^{\overline c}
+\sum_{a\ne \overline c}\text{coef}_X(w_\epsilon(X),a)X^a.
$$
Since noc$(F,G_v)=c$, we get
$$
\text{ord}_X F(X^\nu,y(X))=(\nu/M_v)\text{int}(F,G_v).
$$
Therefore ord$_X F(X^{\nu},y(X)))<0$.
It follows that we are in the situation of (9.33) with $c=c'$.
Consequently by (9.33.2) we see that ${\rm ord}_XG(X^\nu,y(X))<0$ and
\begin{equation*}
\frac{{\rm ord}_X F(X^\nu,y(X))}{{\rm ord}_X G(X^\nu,y(X))}
=\frac{N}{M}
=\frac{\sum_{j\in J}\widehat d(F_j,c)}
{\sum_{l\in J'}\widehat d(G_l,c)}.
\tag{1}
\end{equation*}
Clearly $\epsilon\in I'$ and for all $i\in I$ we have
${\rm coef}_X(z_i(X),c\nu)\ne{\rm coef}_X(w_\epsilon(X),c\nu)$.
Therefore by (9.33.3) we see that
$(F_j)_{j\in J},(G_l)_{l\in J'}$ is a saturated equilateral bisequence of 
$(F,G)$ with diameter $c$. In particular, 
for all $i\in I$ and $e\in I'$ we have
$$
y(X)=\lambda X^{\overline c}
+\sum_{a\ne \overline c}\text{coef}_X(z_i(X),a)X^a
$$
and
$$
y(X)=\lambda X^{\overline c}
+\sum_{a\ne \overline c}\text{coef}_X(w_e(X),a)X^a.
$$
Likewise, for all $j\in J$ and $l\in J'$ we have
noc$(F,G_l)=c$ and noc$(F_j,G)=c$, and hence
\begin{equation*}
\begin{cases}
\text{ord}_X F(X^\nu,y(X))=(\nu/M_l)\text{int}(F,G_l)\\
\text{and}\\
\text{ord}_X G(X^\nu,y(X))=(\nu/N_j)\text{int}(F_j,G).
\end{cases}
\tag{2}
\end{equation*}
By (1) and (2) we get (9.32.14) and (9.33.15), and hence our bisequence
is balanced.

\centerline{}

{\bf Equicognate Lemma (9.35).}
{\it Let $F$ and $G$ in $R$ be $Y$-monic of $Y$-degrees $N>0$ and $M>0$ 
respectively. Then any balanced equicognate bisequence of $(F,G)$ is 
well-balanced.}

\centerline{}

PROOF. Let $(F_j)_{j\in J},(G_l)_{l\in J'}$ be a balanced equicognate 
bisequence of $(F,G)$ with diameter $c$.
By (9.32.8) to (9.32.10) there exists a unique positive integer $B$ such that
for all $j\in J$ and $l\in J'$ we have $N_j=B\widehat d(F_j,c)$ and
$M_l=B\widehat d(G_l,c)$. Hence by (9.32.15) we get
$$
\frac{\sum_{j\in J}N_j}{\sum_{l\in J'}M_l}
=\frac{N}{M}.
$$
By (9.32.14) there exist unique negative rational numbers $N'$ and $M'$
with $MN'=NM'$ such that for all $j\in J$ and $l\in J'$ we have
$N_j=N'\text{int}(F_j,G)$ and $M_l=M'\text{int}(F,G_l)$. 
Namely, to obtain the existence of $M'$, fixing some $j\in J$ and letting
$M'=\frac{MN_j}{N\text{int}(F_j,G)}$, by (9.32.14) we see that for all
$l\in J'$ we have $\frac{M_l}{\text{int}(F,G_l)}=M'$. By symmetry we get the
existence of $N'$. The uniqueness follows in a similar manner.
This proves (9.32.18). From this and the above display we see that
$$
\sum_{j\in J}\text{int}(F_j,G)=\sum_{l\in J'}\text{int}(F,G_l).
$$
In view of (9.32.8) to (9.32.10),
upon letting $E=\text{max}((N_j)_{j\in J},(M_l)_{l\in J'})$
and $D=\text{min}(N_j)_{j\in J}$ and $D'=\text{min}(M_l)_{l\in J'}$,
we have that: 

(i) if $D\ne E$ then there is a unique $s\in J$ such that
$N_s|E$ with $N_s<E=N_j=M_l$ for all 
$j\in J\setminus\{s\}$ and $l\in J'$, and

(ii) if $D'\ne E$ then there is a unique $s'\in J'$ such that
$M_{s'}|E$ with $M_{s'}<E=N_j=M_l$ for all 
$j\in J$ and $l\in J'\setminus\{s'\}$.

Therefore the bisequence is well-balanced.

\centerline{}

{\bf Equilateral Lemma (9.36).}
{\it Every equilateral sequence in $R^{\natural}$ is equicognate.}

\centerline{}

PROOF. Let $(f_i)_{1\le i\le r}$ be an equilateral sequence in
$R^{\natural}$ of diameter $c$, and let $n_i$ be the $Y$-degree of $f_i$. 
By (9.32.0) we see that
$\widehat m(f_i)/n_i\le c$ for $1\le i\le r$.
If $\widehat m(f_i)/n_i=c$ for $1\le i\le r$
then clearly the sequence is cognate. 
Assuming the contrary, after suitable relabelling 
we may henceforth suppose that $\widehat m(f_1)/n_1<c$. 
If $\widehat m(f_i)/n_i<c$ for $2\le i\le r$ then clearly the sequence is 
overcognate, and if $\widehat m(f_i)/n_i=c$ for $2\le i\le r$ 
then clearly the sequence is subcognate.
Assuming the contrary, we must have $r\ge 3$ and after suitable relabelling 
we may henceforth also suppose that 
$\widehat m(f_2)/n_2<c=\widehat m(f_3)/n_3$. 
Now both $f_1$ and $f_2$ are subconjugates of $f_3$, and hence by applying
the coefficient equation of (9.32.7) to the pairs $(f_1,f_3)$ and
$(f_2,f_3)$ we get noc$(f_1,f_2)>c$. This contradicts the fact that $c$ is
the diameter of our equilateral sequence. Thus our sequence is
equicognate.

\centerline{}

{\bf First Packet Lemma (9.37).}
{\it Let $F$ and $G$ in $R$ be $Y$-monic of $Y$-degrees $N>0$ and $M>0$ 
respectively, and assume that $J(F,G)\in k((X))$. 
Then we have the following.

(9.37.1) If $l^*$ is an integer with $1\le l^*\le\chi(G)$ 
and {\rm int}$(F,G_{l^*})<0$, then upon letting
$$
\begin{cases}
c={\rm noc}(F,G_{l^*}),\\
J=\{j:1\le j\le\chi(F)\text{ with }{\rm noc}(F_j,G_{l^*})\ge c\},\\
J'=\{j:1\le l\le\chi(G)\text{ with }{\rm noc}(G_l,G_{l^*})\ge c\},
\end{cases}
$$
we have that $c$ is a rational number such that $c\nu$ is an integer and 
$(F_j)_{j\in J},(G_l)_{l\in J'}$ is a packet of $(F,G)$ with diameter $c$,
and hence in particular $J\ne\emptyset$ and $l^*\in J'$.

(9.37.2) If $j^*$ is an integer with $1\le j^*\le\chi(F)$ 
and {\rm int}$(F_{j^*},G)<0$, then upon letting
$$
\begin{cases}
c={\rm noc}(F_{j^*},G),\\
J=\{j:1\le j\le\chi(F)\text{ with }{\rm noc}(F_{j^*},F_j)\ge c\},\\
J'=\{j:1\le l\le\chi(G)\text{ with }{\rm noc}(F_{j^*},G_l)\ge c\},
\end{cases}
$$
we have that $c$ is a rational number such that $c\nu$ is an integer and 
$(F_j)_{j\in J},(G_l)_{l\in J'}$ is a packet of $(F,G)$ with diameter $c$,
and hence in particular $j^*\in J$ and $J'\ne\emptyset$.}

\centerline{}

PROOF. (9.37.1) follows from (9.34) to (9.36). By symmetry, (9.37.2)
follows from (9.37.1).

\centerline{}

{\bf Second Packet Lemma (9.38).}
{\it Let $F$ and $G$ in $R$ be $Y$-monic of $Y$-degrees $N>0$ and $M>0$ 
respectively, and assume that $J(F,G)\in k((X))$. 
Let $\overline J=\{j:1\le j\le\chi(F)\text{ with }{\rm int}(F_j,G)<0\}$
and $\overline J'=\{j:1\le l\le\chi(G)\text{ with }{\rm int}(F,G_l)<0\}$.
Then we have the following.

(9.38.1) There exists a nonnegative integer $r$ together with disjoint 
partitions $\overline J=\coprod_{1\le i\le r}J(i)$ and
$\overline J'=\coprod_{1\le i\le r}J'(i)$ of $\overline J$ and
$\overline J'$ into pairwise disjoint nonempty subsets such that
$(F_j)_{j\in J(i)},(G_l)_{l\in J'(i)}$ is a packet of $(F,G)$ for
$1\le i\le r$. Moreover these partitions are unique up to order.

(9.38.2) If ${\rm int}(F,G)\ne 0$ and for $1\le j\le\chi(F)$ 
and $1\le l\le\chi(G)$ we have 
${\rm int}(F_j,G)\le 0$ and ${\rm int}(F,G_l)\le 0$,
then $\overline J\ne\emptyset\ne\overline J'$ and in the notation of
(9.38.1) we have $r>0$.}

\centerline{}

PROOF. (9.38.1) follows from (9.37). (9.38.2) is obvious.

\centerline{}

{\bf Remark (9.39).} Let $F$ and $G$ in $R$ be $Y$-monic of 
$Y$-degrees $N>0$ and $M>0$ respectively. 
To elucidate the hypotesis of (9.38.2), and in analogy with the notion of
minint introduced in Section 1,
we define the {\bf strict minimal intersection} of $F$ and $G$ by putting
$$
\text{sminint}(F,G)=\text{min}_{(u,v)\in k^2}\text{int}(F-u,G-v).
$$
Now if int$(F,G)={\rm sminint}(F,G)\ne 0$ then clearly the hypothesis of
(9.38.2) is satisfied, i.e.,
${\rm int}(F,G)\ne 0$ and for $1\le j\le\chi(F)$ 
and $1\le l\le\chi(G)$ we have 
${\rm int}(F_j,G)\le 0$ and ${\rm int}(F,G_l)\le 0$,
Continuing the discussion of Remark (2.4), let us
say that the pair $(F,G)$ is {\bf generic} to mean that
int$(F,G)=\text{sminint}(F,G)$.
By taking indeterminates $U,V$ over $R$, we can consider intersection
multiplicities in $k^*((X))[Y]$ where $k^*$ is an algebraic closure of 
$k(U,V)$. Then, assuming GCD$(F,G)=1$, we get
$$
\text{Res}_Y(F-U,G-V)=\Theta(U,V)X^i+\text{terms of $X$-degree}>i
$$
where 
$$
i=\text{int}(F-U,G-V)\quad\text{ and }\quad
0\ne\Theta(U,V)\in k[U,V].
$$
It follows that 
$$
\text{int}(F-U,G-V)=\text{sminint}(F,G)
$$
and hence for any $(u,v)\in k^2$ we have:
$$
(F-u,G-v)\text{ is generic }\Leftrightarrow\Theta(u,v)\ne 0.
$$
It can be shown that if $(F,G)\sim_m(f,g)$ with $f,g$ in $k[X,Y]$ then
the field degree $[k(X,Y);k(f,g)]$ equals the intersection multiplicity
int$(F-U,G-V)$. Moreover,
if $f,g$ is an automorphic pair, i.e., if $k[f,g]=k[X,Y]$, then by \cite{Ab2} 
we see that $F,G$ are irreducible over 
$k((X))$ and their roots coincide upto the last characteristic term, i.e.,
they are pseudocognates of each other in the sense of (9.32).
This motivates (9.38) where we showed that, under certain condition, 
$J(F,G)\in k((X))$ implies that most branches of $F$ and $G$ can be 
partitioned into packets $(F_1,\dots,F_r,G_1,\dots,G_s)$ whose members are 
pseudocognates of each other. Thus (9.38) may be viewed as a small
contribution to the jacobian problem.

\centerline{}

\centerline{\bf Section 10: Enhanced Newton Polygon}

To simplify the statement of Main Lemma (9.31), in the notation of Section 9,
assuming $N>0$, we let 
the ENP$(F)=$ the {\bf Enhanced Newton Polygon} of $F$ to consist of
the two sequences
$$
(O_i(F))_{1\le 0\le\iota(F)}
\quad\text{ and }\quad
(P_i^{(F)})_{1\le i\le\iota(F)}.
$$
By (9.5) it follows that 
\begin{equation*}
\text{ENP$(F)$ determines UNP$(F)$.}
\tag{10.1}
\end{equation*}
Assuming $N>0$ and $M>0$,
for $0\le j\le\text{min}(\iota(F),\iota(G))$ we say that 
ENP$(F)$ and ENP$(G)$ are {\bf $j$-step parallel}, 
in symbols we write $\text{ENP}(F)||_j\text{ENP}(G)$, if
$$
\begin{cases}
MO_0(F)=NO_0(G),\\
O_i(F)=O_i(G)\text{ for }1\le i\le j,\text{ and}\\
\text{$P_i^{(F)}$ and $P_i^{(G)}$ are related for $1\le i<j$}.
\end{cases}
$$
Moreover, we say that ENP$(F)$ and ENP$(G)$ are {\bf parallel}, 
in symbols we write $\text{ENP}(F)||\text{ENP}(G)$, if
$$
\iota(F)=\iota(G)\text{ and }
\text{ENP}(F)||_{\iota(F)}\text{ENP}(G).
$$
Likewise, we say that ENP$(F)$ is {\bf smaller} than ENP$(G)$, 
in symbols we write $\text{ENP}(F)<\text{ENP}(G)$, if
$$
\begin{cases}
\widehat O(F)<\widehat O(G)\text{ with deg}_Y\widehat P^{(G)}=1,\text{ and}\\
\text{either }\iota(F)=\iota(G)\text{ with } 
\text{ENP}(F)||_{\iota(F)-1}\text{ENP}(G)\\
\qquad\qquad\qquad\qquad\text{ and }
M\text{deg}_Y\widehat P^{(F)}=N\text{deg}_Y\widehat P^{(G)},\\
\,\;\quad\text{or }\iota(F)=\iota(G)-1\text{ with }
\text{ENP}(F)||_{\iota(F)}\text{ENP}(G).
\end{cases}
$$
Finally, we say that ENP$(F)$ and ENP$(G)$ are {\bf pseudoparallel},
in symbols we write $\text{ENP}(F)|.|\text{ENP}(G)$, if
either $\text{ENP}(F)||\text{ENP}(G)$ or $\text{ENP}(F)<\text{ENP}(G)$ or
$\text{ENP}(G)<\text{ENP}(F)$. 
By (9.5) it follows that 
\begin{equation*}
\begin{cases}
\text{for }0\le j\le\text{min}(\iota(F),\iota(G))\text{ we have:}\\
\text{ENP}(F)||_j\text{ENP}(G)\\
\Leftrightarrow\text{UNP}(F)||_j\text{UNP}(G)\text{ and}\\
\quad\text{$P_i^{(F)}$ and $P_i^{(G)}$ are related for $1\le i<j$}
\end{cases}
\tag{10.2}
\end{equation*}
and
\begin{equation*}
\begin{cases}
\text{ENP}(F)||\text{ENP}(G)\\
\Leftrightarrow\text{UNP}(F)||\text{UNP}(G)\text{ and}\\
\quad\text{$P_i^{(F)}$ and $P_i^{(G)}$ are related for $1\le i<\iota(F)$}
\end{cases}
\tag{10.3}
\end{equation*}
and
\begin{equation*}
\begin{cases}
\text{ENP}(F)<\text{ENP}(G)\\
\Leftrightarrow\text{UNP}(F)<\text{UNP}(G)\text{ and}\\
\quad\text{$P_i^{(F)}$ and $P_i^{(G)}$ are related for $1\le i<\iota(G)-1$.}
\end{cases}
\tag{10.4}
\end{equation*}
Consequently, (9.31) may be restated in the following equivalent form:

\centerline{}

{\bf Main Proposition (10.5).} {\it Let $F$ and $G$ in $R$ be 
of $Y$-degrees $N>0$ and $M>0$ respectively,
and assume that $J(F,G)\in k((X))$ and 
$MO_0(F)=NO_0(G)$. Also assume that
either {\rm int}$(F,Y)<0$ or {\rm int}$(G,Y)<0$. Then
${\rm ENP}(F)|.|{\rm ENP}(G)$.  Moreover, if
$\widehat O(F)=\widehat O(G)$ then ${\rm ENP}(F)||{\rm ENP}(G)$
and hence in particular $(M){\rm int}(F,Y)$ $=$ $(N){\rm int}(G,Y)$
and $M\widehat L(F)=N\widehat L(G)$.}

\centerline{}

{\bf Remark (10.6).}
Let us define the {\bf average} or {\bf postfinal root order} of $F$ by
putting
$$
\widetilde O(F)=(1/N)\text{ord}_X F(X,0).
$$
Now the equation
$(M){\rm int}(F,Y)=(N){\rm int}(G,Y)$ in (9.10), (9.31) and (10,5) may be 
restated as saying $\widetilde O(F)=\widetilde O(G)$.
We may postaugument the $O$-sequence by declaring that
$O_{\iota(F)+1}(F)=\widetilde O(F)$ and noting that then:
$\text{ENP}(F)||\text{ENP}(G)\Rightarrow
\text{UNP}(F)||\text{UNP}(G)\Rightarrow\iota(F)=\iota(G)\text{ and }
O_j(F)=O_j(G) \text{ for }0\le j\le\iota(F)+1$.
Finally, we may close-up CNP$(F)$ by its 
$(\iota(F)+1)$-th line whose slope is $\widetilde O(F)$ and which joins the
point $(\Lambda_1(F),L_1(F))$ to the point
point $(\widetilde\Lambda(F),\widetilde L(F))$, with the understanding
that if $\widetilde O(F)=\infty$ then this is the half-infinite horizontal line
emanating from the point $(\Lambda_1(F),L_1(F))$ and going to infinity
on the right. This $(\iota(F)+1)$-th line may be called the {\bf hypotenuse} 
of CNP$(F)$. 

\centerline{}

{\bf Remark (10.7).} The whole game of the Newton Polygon may be redone
by starting with the polynomial
$P^{(F,c)}=P^{(F,c)}(Y)=\text{inco}_X F(X^\nu,YX^{c\nu})\in k[Y]$
and noting that $O_1(F)<\dots<O_{\iota(F)}(F)$ are exactly those rational
numbers $c$ for which this is a true polynomial, i.e., has at least
two terms. Now  put $P_i^{(F)}=P^{(F,O_i)}$ with
$L_i(F)=\text{deg}_Y P_i^{(F)}$ and note that
$L_{i+1}(F)=\text{ord}_Y P_i^{(F)}$. Special adjustments have to be made
if $F(X,0)=0$.

\centerline{}

\centerline{\bf Section 11: Concordance with Homology Rank}
\centerline{\bf and the Numbers of Milnor and Tjurina}

\centerline{}

As we said in Section 3, the proof of Formula (3.3.2) follows from 
Dedekind's Theorem which says that (the ideal generated by) the derivative 
equals the conductor times the different; see page 65 of \cite{Ab4} where 
it is paraphrased in the geometric aphorism: 
the discriminant locus is the union of
the branch locus and the projection of the singular locus. 
Identity (4.8) may be thought of as a modified version of this, 
and may be codified in the algebraic aphorism: the 
affine derived size equals the modified affine conductor size plus the 
modified affine different size plus the degree minus one.  
Thus for any $f\in R_2$ we call 
int$(f,f_Y;\mathcal A)$ and int$(f_X,f_Y;\mathcal A)$
the {\bf affine derived size} of $f$ and
the {\bf modified affine conductor size} of $f$ respectively; in the 
algebraic aphorism we are calling $\overline\beta(f;\mathcal A)$ the
{\bf modified affine different size} of $f$. As abbreviations we put
$$
\epsilon(f;\mathcal A)=\text{int}(f,f_Y;\mathcal A)
\quad\text{ and }\quad
\mu(f;\mathcal A)=\text{int}(f_X,f_Y;\mathcal A)
$$
where these are nonnegative integers or infinity.
Now (4.8) says that if $f\in R_2$ is $Y$-monic of $Y$-degree $N>0$ with
gcd$(f_Y,f-c;\mathcal A)=1$ for all $c\in k$ then
\begin{equation*}
\epsilon(f;\mathcal A)=\mu(f;\mathcal A)+\overline\beta(f;\mathcal A)+(N-1)
\tag{11.1}
\end{equation*}
where all the terms are nonnegative integers.
By analogy, for any $F\in R$ we put
$$
\epsilon(F)=\text{int}(F,F_Y)
\quad\text{ and }\quad
\mu(F)=\text{int}(F_X,F_Y)
$$
where these are integers or infinity and we call them the {\bf derived size} 
of $F$ and the {\bf modified conductor size} of $F$ respectively. 
For any $f\in R_2$ we put
$$
\mu_0(f;\mathcal A)
=\sum_{\{Q\in\mathcal A:f_Q(0,0)=0\}}\text{int}(f_X,f_Y;Q)
$$
and
$$
\overline\mu(f;\mathcal A)
=\sum_{\{Q\in\mathcal A:f_Q(0,0)\ne0\}}\text{int}(f_X,f_Y;Q)
$$
and call these the {\bf restricted conductor size} of $f$ and the
{\bf corestricted conductor size} of $f$ and we note that then
$$
\mu(f;\mathcal A)=\mu_0(f;\mathcal A)+\overline\mu(f;\mathcal A)
=\sum_{\lambda\in k}\mu_0(f-\lambda;\mathcal A)
$$
where all these quantities are nonnegative integers or infinity.
For any $f\in R_2$ we put
$$
\rho(f)=\overline\mu(f;\mathcal A)+\overline\beta(f;\mathcal A)
$$
and we call $\rho(f)$ the {\bf rank} of $f$ and note that it is a
nonnegative integer or infinity.
Now a paraphrase of (11.1) says that
if $f$ is $Y$-monic of $Y$-degree $N>0$ with
gcd$(f_Y,f-c;\mathcal A)=1$ for all $c\in k$ then
\begin{equation*}
\epsilon(f;\mathcal A)=\mu_0(f)+\rho(f)+(N-1)
\tag{11.2}
\end{equation*}
where all the terms are nonnegative integers.
For any $F\in R$ we put
$$
\overline\chi(F)=\chi(F)-1
$$
and call this the {\bf decreased branch number} of $F$, and we note that it
is an integer $\ge -1$. For any $f\in R_2$ we put
$$
\overline\chi(f;\mathcal A)
=\sum_{\{Q\in\mathcal A:f_Q(0,0)=0\}}\overline\chi(f_Q)
$$
and
$$
\overline\chi(f;\mathcal P)
=\sum_{\{Q\in\mathcal P:f_Q(0,0)=0\}}\overline\chi(f_Q)
$$
and call these the {\bf decreased affine branch number} of $f$ and
the {\bf decreased projective branch number} of $f$ respectively, and note
that they are nonnegative integers or infinity.
If $f\in R_2$ is $Y$-monic of $Y$-degree $N>0$ with
gcd$(f_Y,f-c;\mathcal A)=1$ for all $c\in k$ then by (3.3.3) we see that
\begin{equation*}
\mu_0(f;\mathcal A)+\overline\chi(f;\mathcal A)=2\delta(f;\mathcal A)
\tag{11.3}
\end{equation*}
where all the terms are nonnegative integers,
and hence by (4.10) we get
\begin{equation*}
(N-1)(N-2)+\overline\chi(f;\mathcal P)
=2\delta(f;\mathcal P)+\rho(f)
\tag{11.4}
\end{equation*}
where all the terms are nonnegative integers. 
In view of the genus formula (5.3),
by (11.4) we see that 
if $f\in R_2$ is irreducible $Y$-monic of $Y$-degree $N>0$ with
gcd$(f_Y,f-c;\mathcal A)=1$ for all $c\in k$ then
\begin{equation*}
\rho(f)=2\gamma(f)+\overline\chi(f;\mathcal{P}) 
\tag{11.5}
\end{equation*}
and therefore in this situation our rank $\rho(f)$ coincides with 
Abhyankar-Sathaye's rank $r(f)$ introduced in their paper \cite{ASa}.
With the assumptions as in (11.5), as was pointed out in \cite{ASa},
in the complex case, $\rho(f)$ coincides with
the {\bf rank of the first homology group} of $f$, i.e., of the point-set
$\{(u,v)\in\mathbb C^2:f(u,v)=0\}$.

Formula (3.3.3) can be paraphrased by saying that if $F\in R_0=k[[X,Y]]$ is 
$k$-distinguished of $Y$-degree $N>0$ with rad$(F)=F$ then
\begin{equation*}
\mu(F)+\overline\chi(F)=2\delta(F)
\tag{11.6}
\end{equation*}
where all the terms are nonnegative integers.
In the complex case, a topological proof of (11.6) was given by
Milnor \cite{Mil}, and $\mu(F)$ is sometimes called the {\bf Milnor number}
of $F$.

Continuing with $F\in R_0$ which is $k$-distinguished of 
$Y$-degree $N>0$ with rad$(F)=F$,
and recalling that $B(F)=R_0/(FR_0)$,
we define the nonnegative integer $\tau(F)$ by putting
$$
\tau(F)=\text{the length of the ideal in $B(F)$ generated by the images
of $F_X$ and $F_Y$}
$$
and we call this the {\bf torsion size} of $F$.
It is easily seen that $\mu(F)$ is the length of the ideal in $R_0$
generated by $F_X$ and $F_Y$, and hence $\mu(F)\ge\tau(F)$.
In \cite{Zar} it is shown that if $F$ is irreducible in $R_0$ then
$\tau(F)$ is the length of the torsion submodule of the module of differentials
of $B(F)$. In that paper, Zariski gives an interesting characterization
of those irreducible $F$ for which $\tau(F)=2\delta(F)$.
In the complex case, $\tau(F)$ is sometimes called
the {\bf Tjurina number} of $F$.

\centerline{}

Since in (3.3) we assumed $F$ to be $k$-distinguished, instead of quoting 
(3.3.3) in the above proof of (11.3) we should have really quoted
(11.8.3) proved below. Note that, for any ideal $I$ in $R_0$
the $k$-vector-space dimension of $R_0/I$ is denoted by $[R_0/I:k]$,
and we have the well-known implication:
\begin{align*}
&F,G\text{ in $R_0$ where $F$ is $k$-distinguished}\\ 
&\Rightarrow\text{int}(F,G)=[R_0/(F,G)R_0:k]=\text{int}(G,F).
\tag{11.7}
\end{align*}
Now let us prove the:

\centerline{}

{\bf Suplemented conductor-derivative formula (11.8).} {\it Let $F\in R_0$ be
$k$-distinguished of $Y$-degree $N>0$ with {\rm rad}$(F)=F$. 
Let $H=H(X,Y)\in R_0$ be such that $H=UF$ where $U=U(X,Y)\in R_0$ with
$U(0,0)\ne 0$. Let $V=V(X,Y)\in R_0$ and $W=W(X,Y)\in R_0$ with 
$V(0,0)\ne 0\ne W(0,0)$ be such that
$VH_Y$ is $k$-distinguished of $Y$-degree $N-1$, $WH_X=0$ if $H_X=0$, and if
$H_X\ne 0$ then
$WH_X=X^a[Y^b+c_1(X)Y^{b-1}+\dots+c_b(X)]$
with nonnegative integers $a,b$ and elements $c_1(X),\dots,c_b(X)$ in $k[[X]]$
for which $c_1(0)=\dots=c_b(0)=0$.}
($V$ and $W$ exists by the Weierstrass Preparation Theorem).
{Then
\begin{align*}
{\rm int}(H_X,VH_Y)
&={\rm int}(F_X,F_Y)={\rm int}(F,F_Y)-N+1\\
&={\rm int}(H,VH_Y)-N+1
\tag{11.8.1}
\end{align*}
and
\begin{equation*}
{\rm int}(H,VH_Y)-N+1={\rm int}(F,F_Y)-N+1=2\delta(F)-\chi(F)+1
\tag{11.8.2}
\end{equation*}
and
\begin{equation*}
{\rm int}(WH_X,VH_Y)={\rm int}(F_X,F_Y)=2\delta(F)-\chi(F)+1
\tag{11.8.3}
\end{equation*}
where all the terms in the above three items are integers.}

\centerline{}

PROOF. By taking $(VH_Y,H)$ for $(F,G)$ in (2.2) we see that
$$
\text{int}(VH_Y,H_X)=\text{int}(VH_Y,H)-N+1+\beta(VH_Y,H)
$$
where each term is an integer. By (3.1) we have $\beta(VH_Y,H)=0$, and clearly
$\text{int}(VH_Y,H_X)=\text{int}(H_X,VH_Y)$ and
$\text{int}(VH_Y,H)=\text{int}(H,VH_Y)$, and hence by the above display we get
\begin{equation*}
\text{int}(H_X,VH_Y)=\text{int}(H,VH_Y)-N+1
\tag{1}
\end{equation*}
where each term is an integer. By taking $U=1$ in (1) we get the equation
\begin{equation*}
\text{int}(F_X,F_Y)=\text{int}(F,F_Y)-N+1
\tag{2}
\end{equation*}
where each term is an integer. In view of (11.7), by the derivative formula  
$H_Y=U_YF+UF_Y$ we see that
\begin{equation*}
\text{int}(H,VH_Y)=\text{int}(F,F_Y).
\tag{3}
\end{equation*}
In view of (11.7), by (1), (2), and (3) we get the equations
\begin{equation*}
\text{int}(WH_X,VH_Y)=\text{int}(H_X,VH_Y)=\text{int}(F_X,F_Y)
\tag{4}
\end{equation*}
and the equations (11.8.1).  By Dedekind's Theorem 
(see pages 65 and 150 of \cite{Ab4}) we have
\begin{equation*}
\text{int}(F,F_Y)-N+1=2\delta(F)-\chi(F)+1.
\tag{5}
\end{equation*}
By (3) and (5) we get (11.8.2).
By (2), (4), and (5) we get (11.8.3).


\begin{thebibliography}{999}

\bibitem[Ab1]{Ab1}
S. S. Abhyankar,
{\it Algebraic Space Curves,}
Montreal Lecture Notes, 1971.

\bibitem[Ab2]{Ab2}
S. S. Abhyankar,
{\it On the semigroup of a meromorphic curve, Part I,}
Proceedings of the International Symposium of Algebraic Geometry, Kyoto,
1977, pages 240-414.

\bibitem[Ab3]{Ab3}
S. S. Abhyankar,
{\it Lectures On Expansion Techniques In Algebraic Geometry,}
Tata Institute of Fundamental Research, Bombay, 1977.

\bibitem[Ab4]{Ab4}
S. S. Abhyankar,
{\it Algebraic Geometry For Scientists And Engineers,}
American Mathematical Society, 1990.

\bibitem[Ab5]{Ab5}
S. S. Abhyankar,
{\it Some remarks on the jacobian question,}
Proceedings of the Indian Academy of Sciences
(Mathematical Sciences), vol. 104 (1994), pages 515-542.

\bibitem[AA1]{AA1}
S. S. Abhyankar and A. Assi,
{\it Jacobian of meromorphic curves,}
Proceedings of the Indian Academy of Sciences (Mathematical Sciences),
vol. 109 (1999), pages 117-163.

\bibitem[AA2]{AA2}
S. S. Abhyankar and A. Assi,
{\it Factoring the jacobian,} Contemporary Mathematics,
vol. 266 (2000), pages 1-10.

\bibitem[ASa]{ASa}
S. S. Abhyankar and A. Sathaye,
{\it Uniqueness of plane embeddings of special curves,}
Proceedings of the American Mathematical Society,
vol. 124 (1996), pages 1061-1069.

\bibitem[Ass]{Ass}
A. Assi,
{\it Meromorphic plane curves,} Mathematische Zeitschrift,
vol. 230 (1999), pages 165-183.

\bibitem[Mil]{Mil}
J. Milnor,
{\it Singular Points of Complex Hypersurfaces,}
Princeton University Press, 1968.

\bibitem[Zar]{Zar}
O. Zariski,
{\it Characterization of plane algebroid curves whose module of 
differentials has maximum torsion,}
Proceedings of the National Academy of Sciences,
vol. 56 (1966), pages 781-786.

\end{thebibliography}
\end{document}